
\define\scrO{\Cal O}
\define\Pee{{\Bbb P}}
\define\Zee{{\Bbb Z}}
\define\Cee{{\Bbb C}}
\define\Ar{{\Bbb R}}

\define\dirac{\rlap{/}\partial}

\define\proof{\demo{Proof}}
\define\endproof{\qed\enddemo}
\define\endstatement{\endproclaim}
\define\theorem#1{\proclaim{Theorem #1}}
\define\lemma#1{\proclaim{Lemma #1}}
\define\proposition#1{\proclaim{Proposition #1}}
\define\corollary#1{\proclaim{Corollary #1}}
\define\claim#1{\proclaim{Claim #1}}

\define\section#1{\specialhead #1 \endspecialhead}
\define\ssection#1{\medskip\noindent{\bf #1}}
\NoBlackBoxes
\loadbold

\documentstyle{amsppt}

\topmatter
\title
Algebraic surfaces and Seiberg-Witten invariants
\endtitle
\author {Robert Friedman  and John W. Morgan}
\endauthor
\address Department of Mathematics, Columbia University, New York,
NY 10027, USA\endaddress
\email rf\@math.columbia.edu, jm\@math.columbia.edu  \endemail
\thanks The first author was partially supported by NSF grant
DMS-92-03940. The second author was partially supported by NSF grant
DMS-94-02988.
\endthanks
\endtopmatter


\document
\section{1. Introduction.}

Donaldson theory has shown that there is a deep connection between the
4-manifold topology of a complex surface and its holomorphic geometry
\cite{5}, \cite{6}, \cite{8}, \cite{9}, \cite{14}, \cite{15}. Recently Seiberg
and Witten have introduced a new set of 4-manifold invariants
\cite{33}. These invariants have greatly clarified the structure of 4-manifolds
and have made it possible to prove various conjectures suggested by Donaldson
theory. The invariants take on an especially simple form for K\"ahler
surfaces, as realized by Witten \cite{34} and Kronheimer-Mrowka. For example,
their arguments show that, if $X$ is a minimal
algebraic  surface of general type, then  every orientation-preserving
self-diffeomorphism $f\: X \to X$ satisfies $f^*c_1(K_X) = \pm c_1(K_X)$, where
$K_X$, the canonical line bundle of $X$, is the line bundle whose local
holomorphic sections are holomorphic 2-forms. More generally, in case
$X$ is a K\"ahler surface with
$b_2^+
\geq 3$, the form of the invariants easily establishes some conjectures of
\cite{14},  that the only smoothly embedded 2-spheres of self-intersection $-1$
in $X$ represent the classes of exceptional curves, and that the pullback of
the
canonical class of the minimal model of $X$ is invariant up to sign under every
orientation-preserving diffeomorphism of $X$ (see Section 3 below). In this
paper we shall extend these results to cover the case of
$b_2^+ = 1$. Here, for a K\"ahler surface $X$, $b_2^+(X) = 2p_g(X) +1$, where
$p_g(X)$ is the number of linearly independent holomorphic 2-forms on $X$,
i\.e\.
$\dim H^0(X; K_X)$.  More specifically, we shall prove:

\theorem{1.1} Let $X$ be a minimal surface of general type with $p_g(X) = 0$.
Let $\tilde X$ be a blowup of $X$ at $\ell$ distinct points, and let $E_1,
\dots, E_\ell$ be the exceptional curves on $\tilde X$. Finally let $K_0 \in
H^2(\tilde X; \Zee)$ be the pullback to $\tilde X$ of the canonical class of
the
minimal model $X$ of $\tilde X$. If $f$ is an orientation-preserving
self-diffeomorphism $f\: \tilde X \to \tilde X$, then for
all $i$ there is a $j$ such that $f^*[E_i] = \pm [E_j]$  and
$f^*K_0 =
\pm K_0$.  More generally, let $X$ and $X'$ be two minimal
surfaces of general type satisfying the above hypotheses. Suppose that $\tilde
X$ and $\tilde{X'}$ are blowups of $X$ and $X'$ respectively at distinct
points,
that
$E_1, \dots, E_\ell$ and
$E_1', \dots E_m'$ are the exceptional curves on $\tilde X$ and $\tilde{X'}$
respectively and that $K_0$ and $K_0'$ are the pullbacks to $\tilde X$ and
$\tilde{X'}$  of the canonical classes of $X$ and $X'$. If
$f\: \tilde X\to \tilde{X'}$ is an orientation-preserving diffeomorphism, then
$\ell = m$, for every $i$ there exists a $j$ such that $f^*[E_i'] = \pm [E_j]$
 and $f^*K_0' = \pm K_0$.
\endproclaim

More generally, we can replace embedded 2-spheres of square $-1$ in the above
theorem by more general negative definite 4-manifolds. Here, if $N$ is a
negative definite 4-manifold and $X$ is a K\"ahler surface which is
orientation-preserving diffeomorphic to $M\#N$ for some 4-manifold $M$, then
it is essentially a remark of  Kotschick (see \cite{17}, \cite{20})
that
$N$ has no nontrivial finite covering spaces and in particular $H_1(N; \Zee) =
0$.  Thus, if $N$ is a negative definite summand
of a K\"ahler manifold, $H_1( N; \Zee) = 0$ and $H^2(N; \Zee)$ is torsion free.
Again by a theorem of Donaldson \cite{7}, $H^2(N; \Zee)$ has a basis $\{n_1,
\dots, n_\ell\}$ such that $n_i^2 = -1$ and $n_i\cdot n_j = 0$ if $i\neq j$.
Such a basis is unique up to sign changes and permutation, and we will refer to
the $n_i$ as the {\sl exceptional classes\/} of $N$.

\theorem{1.2} Let $X$ be a minimal surface of general type with  $p_g(X) = 0$.
Let $\tilde X$ be a blowup of $X$ at $\ell$ distinct points, and let $E_1,
\dots, E_\ell$ be the exceptional curves on $\tilde X$. Let $N$ be a closed
oriented negative definite $4$-manifold, and suppose that
$\{n_1,
\dots, n_k\}$ is a basis for $H^2(N; \Zee)$ such that $n_i^2 = -1$ for all $i$
and $n_i\cdot n_j = 0$ if $i\neq j$. If there is an orientation-preserving
diffeomorphism $\tilde X \to M \# N$, then, for every $i$, $n_i = \pm [E_j]$
for some $j$.
\endproclaim

Using the above theorems, the arguments of Witten, Kronheimer, and Mrowka in
case $p_g>0$, and standard material on algebraic surfaces (see e\.g\.
\cite{15}) we can deduce the following corollary:

\corollary{1.3} Let $X$ be a minimal K\"ahler surface which is not rational or
ruled, or equivalently such that the Kodaira dimension of $X$ is at least
zero. Then the conclusions of \rom{(1.1)} and \rom{(1.2)} hold for $X$.
\endproclaim

By the classification of surfaces, a K\"ahler surface $X$ which is not
rational or ruled either satisfies $p_g(X) >0$, $p_g(X) = 0$ and $X$ is
of general type, or $p_g(X) = 0$ and $X$ is elliptic. To establish the
corollary, the case $p_g >0$ is covered by the arguments of
Witten and Kronheimer-Mrowka. The case $p_g=0$ and $X$ of general type is
covered by the above theorems. There remains the case that $p_g=0$ and $X$
is elliptic. The case where $p_g(X) = 0$, $X$ is elliptic, and $b_1(X) = 0$
(essentially the case of Dolgachev surfaces) can be handled by arguments
similar
to the proof of the above theorems.  In the remaining case $p_g(X) = 0$ and
$b_1(X) = 2$. In this paper,  we shall just show by elementary methods that, in
the notation of (1.1), every orientation-preserving self-diffeomorphism
$f$ preserves $\pm K_0$ up to torsion, and similarly for the case of two
different surfaces $X$ and $X'$. In fact, one can also  show that $\pm K_0$
itself is preserved.

Using the above results and the general theory of complex surfaces (not
necessarily K\"ahler), we can easily deduce that the plurigenera $P_n(X)$ of a
complex surface are smooth invariants:

\corollary{1.4} If $X$ and $X'$ are two diffeomorphic complex surfaces, then
$P_n(X) = P_n(X')$ for all $n\geq 1$.
\endproclaim

The method of proof of Theorem 1.1 also yields:

\corollary{1.5} Let $X$ be a K\"ahler surface, not necessarily minimal. If
there
exists a Riemannian metric of positive scalar curvature on $X$, then $X$ is
rational or ruled.
\endproclaim

{\it Acknowledgements\/}: It is a pleasure to thank Stefan Bauer, Ron
Fintushel, Tom Mrowka, Ron Stern, Cliff Taubes, Gang Tian, and Edward Witten
for many stimulating discussions about the new invariants. We are especially
grateful to Peter Kronheimer for showing us how to remove 2-torsion from the
statements of the results by taking differences of Spin${}^c$ structures.

\section{2. Seiberg-Witten invariants for K\"ahler metrics.}

Here we review the general theory of  Seiberg-Witten invariants for K\"ahler
metrics. None of the results in this section are original, and most should
appear in \cite{10}, but we shall sketch some of the proofs for the sake of
clarity. Other expositions of this material and its consequences for K\"ahler
surfaces have appeared in
\cite{4}, \cite{24}. Let
$M$ be a a general closed oriented Riemannian 4-manifold with Riemannian metric
$g$.  First we recall some of the properties of Spin${}^c$
structures on $M$. There is an exact sequence
$$\{1\} \to U(1) \to \operatorname{Spin}^c(4) \to SO(4) \to \{1\},$$
which realizes $\operatorname{Spin}^c(4)$ as a central extension of $SO(4)$. In
particular, by considering the exact sequence of cohomology sets associated
to this central extension, the set of Spin${}^c$ structures on
$M$ lifting the frame bundle (which is nonempty if and only if $w_2(M)$ is
the mod 2 reduction of an integral class) is a principal homogeneous space over
$H^1(M;\Cal U(1))$, where $\Cal U(1) = C^\infty(M)/\Bbb Z$ is the sheaf of
$C^\infty$ functions from $M$ to $U(1)$. Since
$H^1(M;C^\infty(M)/\Bbb Z) \cong H^2(M; \Bbb Z)$, given two
Spin${}^c$  structures $\xi _1, \xi _2$ on
$M$ lifting the frame bundle, their difference $\delta (\xi _1, \xi _2)$ is a
well-defined element of $H^2(M; \Bbb Z)$. In dimension four,
$\operatorname{Spin}^c(4)$ is the subgroup of $U(2) \times U(2)$ consisting
of pairs $(T_1,T_2)$ with $\det T_1 = \det T_2$. Thus there are two natural
homomorphisms $\operatorname{Spin}^c(4) \to U(2)$, and the corresponding
homomorphisms $\operatorname{Spin}^c(4) \to U(1)$ given by taking the
determinant agree. If $\xi$ is a Spin${}^c$ structure on $M$, there are two
associated $U(2)$ bundles $\Bbb S^{\pm}= \Bbb S^{\pm}(\xi)$, and $L = c_1(\Bbb
S^{\pm})$ is a complex line bundle which satisfies $c_1(L) \equiv w_2(M) \bmod
2$. Thus $L$ is characteristic, i\.e\. $c_1(L)$ has mod two reduction equal to
$w_2(M)$. We shall call $L$ the complex line bundle {\sl associated to\/}
$\xi$.  Let
$(L, \xi)$ be a pair consisting of  a characteristic complex line bundle
$L$ on $M$, and a
Spin${}^c$ structure $\xi$ whose associated line bundle is $L$. Of course, the
Spin${}^c$ structure determines the line bundle $L$, but we shall record both
since we shall primarily be interested in  $L$ (and we shall sometimes omit the
$\xi$ from the notation). Conversely, $L$ determines $\xi$ up to 2-torsion in
$H^2(M; \Zee)$, and thus $L$ determines $\xi$ uniquely if there is no
2-torsion in $H^2(M; \Zee)$. In terms of the pairs
$(L_1,
\xi _1)$ and $(L_2, \xi _2)$, we have $2\delta (\xi _1, \xi _2) = c_1(L_1)
-c_1(
L_2)$.  Thus we may define the difference $\dsize\frac{c_1(L_1) -c_1( L_2)}2=
\delta (\xi _1, \xi _2)$,
and this difference is well-defined in integral cohomology (not just modulo
2-torsion).

In case $X$ is a K\"ahler surface, or more generally a 4-manifold with an
almost complex structure, then there is a natural lifting of the reduction
of the structure group of the tangent bundle of $X$ to
$\operatorname{Spin}^c(4)$. Namely, we take the map $U(2) \to U(2) \times
U(2)$ defined by $(\rho _1, \rho _2)$, where $\rho _2 = \operatorname{Id}$
and
$$\rho _1(T) = \pmatrix 1 & 0\\ 0 & \det T \endpmatrix.$$
Thus $\Bbb S^+ = \underline \Cee \oplus K_X^{-1}$ and $\Bbb S^- = T_\Cee$,
where $K_X^{-1}$ is the inverse of the canonical line bundle of the almost
complex structure and $T_\Cee$ is the tangent bundle, viewed as a complex
2-plane bundle. In terms of the bundles of $(p,q)$-forms defined by the almost
complex structure, we may also write this as
$\Bbb S^+ = \Omega ^0_X  \oplus \Omega ^{0,2}_X$ and $\Bbb S^- = \Omega
^{0,1}_X$. For this lift $L = K_X^{-1}$. If we replace $\xi$ by $\xi \otimes
\Xi$, where
$\Xi$ is a $C^\infty$ line bundle on $X$, then we replace $K_X^{-1}$ by $L=
\Xi^{\otimes 2} \otimes K_X^{-1}$.  Thus, for a characteristic complex (but not
necessarily holomorphic) line bundle
$L$, a Spin${}^c$ structure for
$X$ is the same as the choice of a complex line bundle
$\Xi$ with
$\Xi
\otimes \Xi = L \otimes K_X$. In this case, given two such Spin${}^c$
structures,
the difference
$\dsize\frac{c_1(L_1) - c_1(L_2)}2$ is equal to
$c_1(\Xi _1) - c_1(\Xi _2)$. Of course, the obvious choices are $L= K_X^{-1}$
and $L = K_X$, with $\Xi$ the trivial bundle in the first case and $\Xi = K_X$
in the second. We shall refer to these two choices as the {\sl natural\/}
Spin${}^c$ structures.

Recall that the Seiberg-Witten equations associated to a {\rm Spin}${}^c$
structure $\xi$ on a Riemannian manifold $M$ are equations for a pair
$(A,\psi)$ where $A$ is a unitary connection of the determinant line
bundle $L$ of $\xi$ and $\psi$ is a section of the bundle $\Bbb S^+(\xi)$ of
plus
spinors.
The equations are
$$\align
F_A^+&=q(\psi)=\psi\otimes
\psi^*-\frac{|\psi|^2}{2}\operatorname{Id}\\
\dirac_A\psi&=0.
\endalign$$  For a generic
perturbation of these equations the moduli space of gauge equivalence classes
of
solutions forms a compact smooth manifold
${\Cal M}(L,\xi)$ which is oriented by a choice of orientation of
$H^1(M;{\Ar})$ and $H^2_+(M;{\Ar})$.
Furthermore, there is a natural two-dimensional cohomology class
$\mu\in H^2({\Cal M}(L,\xi))$ obtained by dividing the space of solutions
to the Seiberg-Witten equations by the group of based changes of
gauge, and taking the first Chern class of the $S^1$-bundle that arises.

Once we have the Seiberg-Witten moduli space, we can define the Seiberg-Witten
invariant for
$M$, which  is a function
$SW _{M,g}$ which assigns an integer to each pair $(L, \xi)$ as above.
It suffices to evaluate $\mu ^{d/2}$ over the fundamental class of $\Cal
M(L, \xi)$, where $d = \dim \Cal M(L, \xi)$. By definition this integer is zero
if $d$ is odd. Strictly speaking, we also need to choose an orientation for
$H^1(X; \Ar) \oplus H^2_+(M; \Ar)$ to determine the sign of
$SW _{M,g}$, but we shall be a little careless on this point. If
$b_2^+(M) \geq 3$, then
$SW _{M,g}$ is independent of the choice of $g$,
whereas if
$b_2^+(M) = 1$, then $SW _{M,g}$ is only defined for generic $g$ and we
will describe the dependence of
$SW _{M,g}$ on $g$ more precisely later.  If $SW _{M,g}(L, \xi)\neq 0$  we
shall
say that the pair $(L, \xi)$ (or $L$) is a  {\sl basic class} (for $g$).  This
can only happen in case the index of the corresponding linearized equations is
nonnegative. This index is
$$\frac14(L^2 -(2\chi (M) + 3\sigma (M)),$$
where $\chi (M)$ is
the Euler characteristic and $\sigma (M)$ is the signature, and we shall also
refer to this index as the {\sl index\/} of the basic class $L$. If
$M = X$ is a complex surface, then it follows from the Hirzebruch signature
formula that
$2\chi (X) + 3\sigma (X) = K_X^2$. Thus if $L$ is a basic class, then $L^2 \geq
K_X^2$. If all of the basic classes have index zero, then
$M$ is called {\sl of simple type} (for $g$). Finally we note that, even though
the function $SW$ is defined by perturbing the Seiberg-Witten equations, if the
unperturbed Seiberg-Witten equations have no solution for $(L, \xi)$, then $SW
_{M,g}(L, \xi) =0$. If on the other hand the solutions to the unperturbed
Seiberg-Witten equations are transverse in an appropriate sense, then we can
use
the moduli space for the unperturbed equations to calculate $SW _{M,g}(L,
\xi)$.

Recall that a solution $(A, \psi)$ to the  SW equations is {\sl reducible\/}
if
and only if $\psi = 0$ and hence
$F_A^+ = 0$, where $F_A^+$ is the self-dual part of the curvature of a
connection on
$L$. In the case of a K\"ahler metric,  $F_A^+ = 0$ implies that $A$ is a
$(1,1)$-connection (and thus $A$ defines a holomorphic structure on $L$) and
$c_1(L) \cdot \omega = 0$, where $\omega$ is the K\"ahler form. By
the Hodge index theorem, $L^2\leq 0$. Conversely, if
$L$ is a holomorphic line bundle and $c_1(L) \cdot \omega = 0$, then there
exists a $(1,1)$-connection $A$ on $L$ with $F_A^+ = 0$, giving a reducible
solution to the SW equations, and indeed in this case all solutions will be
reducible. We call a basic class $L$ {\sl reducible\/} if all solutions to the
corresponding SW equations are reducible, and {\sl irreducible\/} otherwise. Of
course, reducible solutions to the SW equations for $(L, \xi)$ do not
necessarily imply that $(L, \xi)$ is basic. For a generic K\"ahler metric,  we
can assume that the K\"ahler form
$\omega$ is not orthogonal to any nontorsion class $L\in H^2(X; \Zee)$ with
$L^2
\geq K_X^2$. We shall call a K\"ahler metric whose associated K\"ahler form
$\omega$ satisfies this condition {\sl generic}. Thus for a generic K\"ahler
metric, there exist reducible basic classes  $L$ only if
$c_1(L)$ is zero as an element of $H^2(X; \Ar)$. Of course, these will give
solutions of nonnegative index only if $0 = L^2 \geq K_X^2$.

In the K\"ahler case, we have the following criterion for the SW moduli
space to be nonempty and of nonnegative formal dimension \cite{34},
\cite{10}:

\proposition{2.1} Suppose that $X$ is a K\"ahler surface with K\"ahler form
$\omega$. Then the pairs $(L, \xi)$ of nonnegative index admitting  irreducible
solutions to the SW equations are in one-to-one correspondence with
holomorphic characteristic line bundles $L$, together with a choice of a
holomorphic square root $(K_X\otimes L)^{1/2}$ for the line bundle $K_X\otimes
L$,  satisfying:
\roster
\item"{(i)}"  $L^2 \geq K_X^2$;
\item"{(ii)}" Either $H^0(X; (K_X\otimes L)^{1/2}) \neq 0$ and $\omega \cdot L
< 0$ or $H^0(X; (K_X\otimes L^{-1})^{1/2}) \neq 0$ and $\omega \cdot L
> 0$. \qed
\endroster
\endproclaim

Here the choice of a square root $(K_X\otimes L)^{1/2}$ for $K_X\otimes L$
naturally gives the square root $(K_X\otimes L^{-1})^{1/2} = (K_X\otimes
L)^{1/2} \otimes L^{-1}$ for $K_X\otimes L^{-1}$.

The idea behind the proof of (2.1) is that since $X$ has a complex structure,
$\Bbb S^+(\xi)$ splits as a sum of line bundles $(K_X\otimes L)^{1/2}$ and
$\Omega^{0,2}_X(( K_X\otimes L)^{1/2})$. Thus,
the spinor field $\psi$ decomposes
into components $(\alpha,\beta)$.  The curvature part of the
Seiberg-Witten equations  says that
$$\align
F^{0,2}
&=\bar\alpha\beta\\
(F_A^+)^{1,1} &=\frac{i}{2}(|\alpha|^2-|\beta|^2)\omega,
\endalign$$
where $\omega$ is the K\"ahler form. The Dirac equation for a K\"ahler surface
becomes
$$\bar\partial_A\alpha+\bar\partial^*_A\beta=0.$$
Applying $\bar\partial_A$ to this equation we get
$$\bar\partial_A\bar\partial_A\alpha
+\bar\partial_A\bar\partial^*_A\beta =0.$$
Equivalently,
$$F^{0,2}_A\cdot
\alpha+\bar\partial_A\bar\partial^*_A\beta =0.$$
Since $F^{0,2}_A=\bar\alpha\beta$,  this equation
becomes
$$|\alpha|^2\beta+\bar\partial_A\bar\partial^*_A\beta =0.$$
Taking the $L^2$-inner product with $\beta$  yields
$$\int _X|\alpha|^2|\beta|^2+\|\bar\partial^*_A\beta\|^2_{L^2}=0.$$

It follows that $|\alpha|^2|\beta|^2$ and $\bar\partial^*_A\beta$ are zero.
Thus  $F^{0,2}_A=\bar\alpha\beta =0$. This means that $A$ is a holomorphic
connection and so defines a holomorphic structure on $L$.   Moreover
$\bar\partial^*_A\beta =0$ and so $\bar\partial_A\alpha =0$. Hence $\alpha$ is
a holomorphic section of $(K_X\otimes L)^{1/2}$ and  $\bar\beta$ is a
holomorphic
section of
$(K_X\otimes L^{-1})^{1/2}$.
Since
$\alpha$ and $\bar \beta$ are holomorphic, they do not vanish
on any open subset unless they vanish identically. So either $\alpha = 0$ or
$\beta=0$.  Furthermore,
$$\omega \cdot L = \int _X\omega \wedge \frac{i}{2\pi}F_A^+ = -\frac{1}{4\pi}
\int _X (|\alpha|^2 - |\beta|^2)\omega \wedge \omega,$$
and so
$\alpha$ is not zero if and only if $\omega \cdot L <0$ and
$\beta$ is not zero if and only if $\omega \cdot L >0$.
If $\alpha\neq 0$
then $\alpha$ is a nonzero holomorphic section of $(K_X\otimes L)^{1/2}$.  If
$\beta\neq 0$, then $\bar\beta$ is a nonzero holomorphic section of
$(K_X\otimes L^{-1})^{1/2}$.

We have seen that the conditions listed in Proposition 2.1 are
necessary for a solution.  Let us show that they are sufficient as
well. We shall just consider the case where $\omega \cdot L$ is
negative. The holomorphic structure on $L$ uniquely determines a connection
$A$, once we have chosen a hermitian metric on $L$, by choosing the unique
connection compatible with the holomorphic structure and the metric. Fix an
arbitrary hermitian metric on $L$. Given
$L$ and a nontrivial holomorphic section
$\alpha$ of
$(K_X\otimes L)^{1/2}$, we wish to change the  metric on $L$ until the
curvature part of the Seiberg-Witten equation is satisfied.
If we think of varying the hermitian metric by $\exp\lambda$
for some real valued function $\lambda$, then the equation
we need to solve for $\lambda$ is
$$F_A^++(\bar\partial\partial
\lambda)^+=\frac{i}2 e^\lambda |\alpha|^2\omega.$$ For a $(1,1)$-form $\eta$,
$\eta ^+ = \Lambda \eta \cdot \omega$, where $\Lambda$ is contraction with
$\omega$. We take the pointwise contraction with the K\"ahler form
$\omega$ and obtain
$$\Delta \lambda -\frac{|\alpha|}{4}^2e^\lambda -\frac12*(iF_A^+\wedge
\omega)=0.$$ Here $\Delta$ is the negative definite Laplacian on functions (in
Euclidean space it would be $\sum _i\partial ^2/\partial x_i^2$). According to
results of Kazdan-Warner
\cite{18} (first applied in gauge theory to the vortex equation by Bradlow
\cite{3}), there is a unique solution
$\lambda$ to this equation provided that $\int _XiF_A^+\wedge \omega <0$, which
is just the condition that $\omega \cdot L$ is
negative. Thus, we have seen that for each non-trivial
holomorphic section of
$(K_X\otimes L)^{1/2}$ we can obtain a solution to the Seiberg-Witten
equations with the holomorphic section as the spinor field. This completes the
sketch of the proof of the proposition.

A gauge equivalence between two solutions $(A,\psi)$ and $(A',\psi')$
will be a holomorphic isomorphism between the holomorphic structures
$L$ and $L'$ determined by the two connections. It will also carry the
section $\psi$ to $\psi'$.  Since the only holomorphic automorphisms of
a holomorphic bundle are multiplication by nonzero constant functions, this
implies that under the holomorphic identification $\psi$ and $\psi'$
define the same point of  $\Pee H^0((K_X\otimes L)^{1/2})$. Conversely, it
is easy to check that two sections which agree modulo $\Cee ^*$ define gauge
equivalent solutions.  Thus we may identify the moduli space to the
unperturbed equations with
$$\bigcup _L\Pee H^0((K_X\otimes L)^{1/2}),$$
where we think of $L$ as ranging over all holomorphic structures on a fixed
$C^\infty$ line bundle (the set of all such structures is isomorphic to the
complex torus
$\operatorname{Pic}^0X$). The moduli space may thus be identified with an
appropriate component of the Hilbert scheme of curves on $X$. Of course, there
will be another moduli space corresponding to
$L^{-1}$ as well.

We now divide the study of the
basic classes into two cases: the case where $X$ is minimal and the case where
$X$ is not minimal.

\medskip
{\bf The case of a minimal surface.} For a K\"ahler surface $X$ which is not
rational or ruled,
$X$ is minimal if and only if
$K_X$ is {\sl nef\/}: in other words, for every holomorphic curve $C$ on $X$,
$K_X\cdot C
\geq 0$. If $K_X$ is nef, then  $K_X^2 \geq 0$. The case where $K_X$ is nef and
$K_X^2 >0$ is the case where $X$ is of general type. In this case,
since
$K_X^2 >0$, there are no reducible basic classes for any K\"ahler metric.

\proposition{2.2} With notation as above, suppose  that $X$ is a minimal
surface of general type, i\.e\. suppose that $K_X$ is nef and that $K_X^2 >0$,
and that
$\omega$ is a generic K\"ahler metric. Then the only pairs
$(L, \xi)$ satisfying the conditions
\rom{(i)} and \rom{(ii)} are $L = K_X^{\pm 1}$, with the natural
\rom{Spin${}^c$} structures, i\.e\. the ones corresponding to the square root
$K_X$ of
$K_X\otimes K_X$ and the square root $0$ of $K_X \otimes K_X^{-1}$.
\endproclaim
\proof For the proof we use additive notation for holomorphic line bundles
(which
we could identify with divisor classes on $X$). After replacing
$L$ by $-L$, we may assume that $\omega \cdot L < 0$ and that $K_X+L$ is
effective. We have $(K_X+L) \cdot \omega \geq 0$, $L\cdot \omega < 0$, so
there is an $a\geq 1$ such that $(K_X+aL)\cdot \omega = 0$. By the Hodge index
theorem $(K_X+aL)^2 \leq 0$, with equality only if $K_X + aL$ is numerically
trivial. Thus
$$K_X^2 + 2a(K_X\cdot L) + a^2L^2 \leq 0.$$
On the other hand $K_X$ is nef, so that $(K_X+L)\cdot K_X \geq 0$.
Putting together
$$\align
2aK_X^2 + 2aK_X\cdot L &\geq 0;\\
K_X^2 + 2a(K_X\cdot L) + a^2L^2 &\leq 0,
\endalign$$
we obtain
$$(1-2a)K_X^2 + a^2L^2 \leq 0,$$
or in other words
$$L^2 \leq \frac{2a -1}{a^2}K_X^2.$$
But
$$\frac{2a -1}{a^2} = 1 - \fracwithdelims(){a-1}{a}^2 = 1 - \left(1 -
\frac1a\right)^2,$$
which is decreasing for $a\geq 1$. Thus
$$L^2 \leq \frac{2a -1}{a^2}K_X^2\leq K_X^2.$$
Since $L^2\geq K_X^2$, we have $L^2 = K_X^2$.
For equality to hold we must have $a=1$ and $K_X+L$ must be numerically
trivial.
In this case, $K_X+L$ has a section since $\dsize\frac{K_X+L}2$ has a section.
Thus
$K_X+L$ is the trivial divisor, so $L= -K_X$. Moreover $\dsize\frac{K_X+L}2$
is numerically trivial and has a section as well, so that it is trivial. Thus
the Spin${}^c$ structure corresponds to taking the trivial square root of
$K_X+L=0$.
\endproof

Essentially the same argument shows:

\proposition{2.3} With notation as above, suppose that $K_X$ is nef and
that $K_X^2 =0$. If $L$ is a line bundle satisfying the conditions
\rom{(i)} and \rom{(ii)}, then there exists a rational number $r\leq 1$ such
that $L$ is numerically equivalent to $\pm  rK_X$. Moreover, in case $r= \pm
1$, then in fact $L =\pm K_X$, and the {\rm Spin}${}^c$ structures are again
the natural ones.
\qed
\endproclaim

Note that in all cases we have $L^2 = K_X^2$, in other words $X$ is of simple
type for $g$ if $X$ is minimal and $K_X$ is nef. Of course, so far we have not
actually shown that there are any basic classes. But  in case
$L = \pm K_X$, the value of $SW$ is $\pm 1$. To determine the exact sign, we
need to make a choice of orientation for the moduli space.
The orientation convention we shall follow is this: To
orient the relevant determinant line bundle for a general $4$-manifold $M$, we
must choose an orientation for the vector space $H^1(M; \Ar) \oplus H^2_+(M;
\Ar)
\oplus H^0(M;
\Ar)$, by choosing orientations on  $H^1(M; \Ar)$ and $H^2_+(M; \Ar)$ and using
the standard orientation on $H^0(M; \Ar)$. For a K\"ahler surface $X$,
$$ H^2_+(X; \Ar) \cong
\Ar \cdot \omega \oplus (H^{2,0}(X) \oplus H^{0,2}(X))_\Ar$$ and $H^1(X;
\Ar)\cong (H^{1,0}(X) \oplus H^{0,1}(X))_\Ar$. We choose the orientation
given by taking the standard orientation on $\Ar\cdot \omega$ and
using the isomorphism
$(H^{i,0}(X) \oplus H^{0,i}(X))_\Ar \to H^{0,i}(X)$ to transfer the
usual complex orientation on $ H^{0,i}(X)$ to $(H^{i,0}(X) \oplus
H^{0,i}(X))_\Ar$. We then
have the following result, which follows easily by considering the
linearization
of the SW equations
\cite{10}:

\proposition{2.4} For an arbitrary complex surface $X$, if $g$ is a K\"ahler
metric on
$X$ with K\"ahler form
$\omega$ and $\omega \cdot K_X >0$, then the value of
$SW _{X,g}$ on $-K_X$ for the natural {\rm Spin}${}^c$ structure is $1$ and
$SW_{X,g}(K_X)= (-1)^{q+ p_g}$ for the natural {\rm Spin}${}^c$ structure.
\endproclaim

\demo{Sketch of proof}
Let us consider the elliptic complex associated to the unique solution
for $L=-K_X$. The kernel of the Dirac operator is isomorphic to
$H^0(X;{\Cee})\oplus H^{0,2}(X)$.  The cokernel of the Dirac operator
is
$H^{0,1}(X)$. Of course, $H^2_+(X;i{\Ar})$ is, as an oriented vector space,
isomorphic to $(i\Ar)\cdot \omega\oplus H^{0,2}(X)$.
Clifford multiplication by the solution $(\alpha,0)\in H^0\oplus
H^{0,2}$ induces an orientation-preserving isomorphism $H^1(X;i{\Ar})\to
H^{0,1}(X)$. If $Dq$ is the differential of the quadratic mapping $q$,
then $-Dq$  induces an orientation-preserving
mapping
$$H^{0,2}(X)\to H^{0,2}(X),$$ namely multiplication by
$-\bar\alpha$.
The map $Dq$ also induces a map
$$H^0(X;{\Cee})\to (i\Ar)\cdot \omega$$
which sends $\eta$ to $-i\operatorname{Re}\langle \alpha,\eta\rangle$.
Of course, the action of the stabilizer $S^1$ on the spin fields is by
the opposite of the complex orientation.  Hence, the solution to the
equations modulo the action of the group of changes of gauge
is a single point.  The equations are transverse at this
point and the orientation is plus one.

A similar computation shows that
$$SW_{X,g}(K_X) = (-1)^{q(X)+p_g(X)}.$$
In fact, for any smooth 4-manifold
$M$ and any Spin${}^c$ structure $\xi$ there is a naturally defined opposite
Spin${}^c$ structure
$-\xi$ and we have
$$SW_{M,g}(\xi)=(-1)^{(1-b_1(M)+b_2^+(M))/2}SW_{M,g}(-\xi).$$
\enddemo

To end our discussion of minimal K\"ahler surfaces, we consider the example of
elliptic surfaces. For simplicity, and because it is the most interesting case,
we shall just consider the case of simply connected elliptic surfaces, so that
linear, numerical, and homological equivalence are the same. In this case the
basic classes have a certain multiplicity which need not be one, but which we
shall not compute here.  Let
$X$ be a simply connected elliptic surface, with $p_g(X) = p_g$. Then $X$ has
at
most two multiple fibers
$F_1$ and $F_2$, of multiplicities $m_1$ and $m_2$, say. From the canonical
bundle formula
$$K_X = (p_g-1)f + (m_1-1)F_1 + (m_2-1)F_2,$$
where $f$ is the class of a general
fiber:
$f= m_1F_1 = m_2F_2$. Let $D$ be a divisor which is a rational multiple of
$K_X$ and thus of the fiber  $f$, say $D = rf$, and define $\deg D = r$.
For an arbitrary K\"ahler metric $\omega$, normalized so that $\omega \cdot
f=1$, we have $\deg D =
\omega \cdot D$. The basic classes correspond to line bundles
$L$ such that either $K_X + L = 2D$, where $D$ is effective and $\omega \cdot L
\leq 0$, i\.e\. $\dsize 0 \leq \deg D \leq \frac{\deg K_X}2$,
or $K_X - L = 2D$, where $D$ is effective and $\omega \cdot L
\geq 0$, i\.e\. again we have $\dsize 0 \leq \deg D \leq \frac{\deg K_X}2$.
The effective divisor $D$ on $X$ can be written as $af + bF_1 + cF_2$, where
$a\geq 0$, $0\leq b \leq m_1-1$, and $0\leq c \leq m_2 -1$. If $p_g >0$, $\deg
D \leq \deg K_X$ if and only if $a \leq p_g-1$. In this case it is clear that
$D$ is effective if and only if $K_X - D$ is effective. Thus in case $K_X + L =
2D$, where $D$ is effective and $\dsize 0 \leq \deg D \leq \frac{\deg K_X}2$,
we have
$$L = 2D - K_X = K_X - 2(K_X-D) = K_X - 2D',$$
where $D'$ is effective and $\dsize \frac{\deg K_X}2 \leq \deg D' \leq \deg
K_X$. Likewise if $K_X - L = 2D$, where $D$ is effective and $\omega \cdot L
\geq 0$, then $L= K_X - 2D$, where $D$ is effective and $\dsize 0 \leq \deg D
\leq \frac{\deg K_X}2$. Thus we see that in all cases the basic classes are
exactly the classes of the form $K_X - 2D$, where $D$ is effective and $0 \leq
\deg D \leq \deg K_X$. In other words, the basic classes are the classes
$(p_g-1-2a)f + (m_1 - 2b-1)F_1 + (m_2 -2c-1)F_2$ for $0\leq a \leq p_g-1$,
$0\leq b \leq m_1-1$, $0 \leq c\leq m_2-1$. These are exactly the
Kronheimer-Mrowka basic classes \cite{21}, \cite{11}, and one can show that the
appropriate multiplicity, up to sign, to attach to the class $K_X-2D$ with
$D =af + bF_1 + cF_2$, with $a \geq 0$,  $0\leq b \leq m_1-1$, and $0\leq c
\leq
m_2 -1$, is
$\dsize \binom{p_g-1}{a}=\binom{h_++ h_-}{h_+}$, where $h_+ = h^0(D)-1$ and
$h_- =h^0(K_X-D)-1$.

A similar analysis holds for the case of an elliptic surface with $p_g =0$,
where
we must simply analyze the conditions on $L$ directly. For example, we obtain
the SW classes $L$ with nonnegative fiber degree (i\.e\. the $L$ such that $L =
rK_X$ with $r\geq 0$ in rational cohomology) by considering the effective
divisors $D$ such that $L=K_X - 2D$ has nonnegative fiber degree. In this case
$\dsize \frac{K_X-L}2 = D$ and the correct multiplicity, up to sign, is one.

\medskip
{\bf The case of a nonminimal surface.}  It is enough to
consider the case where
$\tilde X
\to X$ is a single blowup.

\proposition{2.5} Let $X$ be a K\"ahler surface which is the blowup of a
surface for which the canonical bundle is nef and let
$g$ be a K\"ahler metric on $X$ with  K\"ahler form $\omega$ such that $\omega$
is not orthogonal to any nontorsion class $L\in H^2(X; \Zee)$ with $L^2 =
K_X^2$.
Let
$\tilde X$ be the blowup of $X$ at a point, with $E$ the exceptional divisor,
and
consider a K\"ahler metric $\tilde g$ on
$\tilde X$ corresponding to the K\"ahler form $\tilde \omega = N\omega -E$, $N
\gg 0$. Then:
\roster
\item"{(i)}" Every basic class on $\tilde X$ for $\tilde g$ is irreducible.
\item"{(ii)}" If  $(\tilde L, \tilde \xi)$ is a basic class on
$\tilde X$ for $\tilde g$, then  either $\tilde L = L\pm E$, where
$L$ is an irreducible basic class on $X$ for $g$, or $\tilde L =L\pm E$, where
$L$ is  a reducible basic class on $X$ and so the image of $L$ is zero in
$H^2(X; \Ar)$. If $\tilde L = L+ E$ is a basic class for $\tilde X$, then the
corresponding {\rm Spin}${}^c$ structure $\tilde \xi$, or in other words the
square root $\tilde \Xi$ of $\tilde L+K_{\tilde X} = L+K_X+ 2E$, is $\Xi + E$,
where
$\Xi$ is the square root of $L+K_X$ corresponding to the {\rm Spin}${}^c$
structure $\xi$ on $X$, and likewise if $\tilde L = L-E$, then $\tilde \Xi =
\Xi$. Moreover, the classes
$L\pm E$ where
$L$ is an irreducible basic class on $X$ for $g$ all give basic classes on
$\tilde X$ for $\tilde g$, and in this case
$$SW_{\tilde X, \tilde g}(L\pm E, \tilde \xi) = \pm SW_{X,g}(L, \xi).$$
\item"{(iii)}" If there is a basic class on $\tilde X$ of the form $L\pm E$,
where $L$ is  a reducible basic class on $X$, then $X$ is minimal. Moreover, if
$K_X$ is torsion, then  $L = \pm K_X$. In this case $\pm K_X \pm E$ are
basic classes on $\tilde X$ and $SW_{\tilde X, \tilde g}(\pm K_X \pm E) = \pm
1$ for the natural {\rm Spin}${}^c$ structures.
\item"{(iv)}" Let $\tilde X$ be a K\"ahler surface which is not rational or
ruled. Let $\tilde g$ be a K\"ahler metric on $\tilde X$ whose K\"ahler form
$\tilde \omega = N\omega - \sum _iE_i$, where the $E_i$ are the exceptional
curves on $\tilde X$, $\omega$ is the K\"ahler form of a generic K\"ahler
metric on the minimal model of $\tilde X$, and $N\gg 0$.  For every basic class
$L$ on
$\tilde X$ for $\tilde g$, we have
$L^2=K_{\tilde X}^2$, so that
$\tilde X$ is of simple type for
$\tilde g$.
\endroster
\endproclaim
\proof First we prove (i). Clearly, for all $N \gg 0$, a K\"ahler metric with
K\"ahler form
$N\omega -E$ is generic. Thus the only possible reducible basic classes are
torsion. But a  basic class is characteristic. If there were a torsion
characteristic element in
$H^2(\tilde X; \Zee)$, then every element of $H^2(\tilde X; \Zee)$ would have
even square. This contradicts the fact that $E^2 =-1$.

Next we consider (ii) and (iii). Suppose that $\tilde L$ is a basic class on
$\tilde X$. Possibly after replacing $\tilde L$ by
$-\tilde L$, we can assume that $\tilde L^2
\geq K_{\tilde X}^2 = K_X^2 -1$, $\dsize\frac{\tilde L + K_{\tilde X}}{2}$ is
effective, and $\tilde \omega \cdot \tilde L <0$. We have $K_{\tilde X} =
K_X+E$, and $\tilde L = L + aE$ for some characteristic $L\in H^2(X; \Zee)$ and
odd integer $a$. As $\dsize\frac{\tilde L + K_{\tilde X}}{2}$ is
effective, $\dsize\frac{ L + K_X}{2}$ is
effective as well.  Since $\tilde
L^2 = L^2 - a^2
\geq K_X^2 -1$, and $a$ is odd, we have $L^2 \geq K_X^2$. As $\tilde \omega
\cdot
\tilde L <0$, we have
$N(\omega \cdot L) + a <0$ for all $N\gg 0$. Thus $\omega \cdot L \leq 0$.
First consider the case where
$\omega \cdot L <0$. Then the line bundle $L$ on $X$ satisfies (i) and (ii) of
(2.1). Conversely, starting with  a $L$ on $X$ satisfying (i) and (ii) of
(2.1), the class
$\tilde L = L\pm E$ will satisfy (i) of (2.1). If say $L\cdot \omega <0$, then
$\tilde L
\cdot \tilde \omega <0$ provided that $N\gg 0$. Finally $\dsize\frac{K_{\tilde
X}
+\tilde L}2 = \frac{K_X+L}2$ or $\dsize \frac{K_X+L}2 +E$, and we have natural
isomorphisms
$$H^0(\tilde X;\frac{K_X+L}2 +E) \cong H^0(\tilde X;\frac{K_X+L}2) \cong
H^0(X;\frac{K_X+L}2),$$
where we have used the notation $H^0(X; D)$ to denote the group of sections of
$\scrO_X(D)$. Thus (ii) of (2.1) is satisfied as well, and we see that the
Spin${}^c$ structures are as claimed. We omit the argument  that $SW_{\tilde X,
\tilde g}(L\pm E) =
\pm SW_{X,g}(L)$, where
$\tilde g$ and
$g$ are appropriate generic K\"ahler metrics on $\tilde X$ and
$X$ respectively. In case $L$ is irreducible, this result is established in the
general case via a general blowup formula in \cite{10}. Note that, in case $L=
K_X$, the main case of interest, $L+E = K_X+E = K_{\tilde X}$, and this case is
covered by (2.4) with $X$ replaced by $\tilde X$. The case $K_X-E$ then follows
by the naturality of the function $SW$, since $K_X-E= R^*(K_X+E)$, where $R\:
\tilde X \to \tilde X$ is the diffeomorphism corresponding to reflection in the
class $E \in H^2(\tilde X; \Zee)$.

Now consider the case where $\omega \cdot L=0$.  By the hypothesis that the
metric is generic, $L$ is  zero in rational cohomology. From $\tilde
L^2 = -a^2
\geq K_{\tilde X}^2 = K_X^2 -1$, we see that $K_X^2 \leq 1-a^2 \leq 0$. Since
$L$ is characteristic,
$X$ must in fact be minimal. Hence $K_X^2
\geq 0$ and so $K_X^2 =0$ and  $a=\pm 1$. In this case  $\tilde L^2 =
K_{\tilde X}^2$, so that all of the $\tilde L$ constructed in this way are of
index zero.  We note that (ii) of (2.1) on $\tilde X$ is satisfied if and only
if $\dsize H^0(X;\frac{K_X+L}2) \neq 0$. Thus $L+K_X$ is effective. If $K_X$ is
torsion, then  $L+K_X$ is also zero in rational cohomology. Thus it is the
trivial divisor, and so $L= -K_X$, and the Spin${}^c$ structure is the trivial
square root of $K_X+L$. It follows from (2.4) that $SW_{\tilde X, \tilde g}
( K_X + E) = \pm 1$, and  arguments as in the irreducible case handle
show that $SW_{\tilde X,
\tilde g}(\pm K_X \pm E) =\pm 1$ also.

Lastly we prove (iv). If $X$ is minimal and $L$ is
irreducible, then (2.2) and (2.3) imply that
$L^2 = K_X^2$.  If $L$ is reducible, then $X$ is minimal and $K_X^2 = 0$.
Since $L$ is torsion $L^2 = K_X^2 =0$ in this case as well. By induction on the
number of blowups, and using the above discussion for reducible
$L$, we may assume that
$\tilde X$ is the blowup of a surface $X$ at one point, where
$L^2=K_X^2$ for all basic classes
$L$ on $X$. We will show that the same is true for $\tilde X$. It suffices to
show that basic classes of the form $L + aE$, where $L$ is an irreducible class
on $X$ and $a$ is an odd integer, have square equal to $K_X^2 -1$. Since
$L^2 = K_X^2$ by induction on the number of blowups, $a^2 \leq 1$. Thus $a =
\pm 1$ and  $\tilde L^2= K_{\tilde X}^2$.
\endproof

We note that using
blowups, we can calculate the invariant in case
the moduli space is singular because the
solutions are all reducible. For example, for a K\"ahler surface $X$ with a
generic metric $g$, suppose $X$ is a minimal surface and
$K_X$ is torsion. In this case the basic classes are exactly $\pm K_X$ with
the natural Spin${}^c$ structure, and the value of $SW_{g,X}$ on these classes
is $\pm 1$.  This is analogous to the use of blowups in Donaldson theory to
define ``unstable" invariants, as in Chapter III, Section 8 of \cite{15}.

\section{3. The case where $\boldkey p_{\boldkey g}$ is nonzero.}

In this section we will describe the proofs of the results corresponding
to Theorems 1.1 and 1.2 in case
$p_g(X) \neq 0$. In this case $SW _{X,g}$ does not depend on the choice of
the metric, and the set of basic classes for $g$, which is independent of the
choice of $g$, is a diffeomorphism invariant of $X$. Let
$X$ be a minimal surface of general type, and let $\tilde X$ be a blowup of
$X$.
We may assume that $\tilde X$ is a blowup of $X$ at distinct points. If
$E_1, \dots, E_\ell$ are the exceptional classes of $\tilde X$, and $K_0$ is
the
pullback to $\tilde X$ of $K_X$, then the basic classes are $\pm K_0 + \sum
_{i=1}^\ell \pm E_i$, with the natural Spin${}^c$ structures. Consider the
subset of all expressions of the form
$\dsize\frac{\tilde L _1 - \tilde L _2}2$, where $\tilde L _1$ and $\tilde L
_2$
are distinct basic classes, where we have used the Spin${}^c$ structures to
define the square roots as integral cohomology classes. Here, if $\tilde L _1 =
\pm K_0 + \sum _{a\in A}E_a + \sum _{a\notin A}(-E_a)$, then the square root
$\tilde \Xi _1$ of $\tilde L _1\otimes K_{\tilde X}$ corresponding to the
choice
of Spin${}^c$ structure is
$$\tilde \Xi _1 = \cases K_0 + \sum _{a\in A}E_a, &\text{if $\tilde L_1 = K_0
+ \sum _{a\in A}E_a + \sum _{a\notin A}(-E_a)$,}\\
 \sum _{a\in A}E_a, &\text{if $\tilde L_1 = -K_0
+ \sum _{a\in A}E_a + \sum _{a\notin A}(-E_a)$.}
\endcases$$
Thus the set of difference classes  consists exactly of the elements
$\pm K_0$, $\pm K_0 + \sum _{i\in A} \pm E_i$, where $A$ is a proper subset of
$\{1, \dots, \ell\}$, or $\sum _{i\in A}\pm E_i$, where $A$ is here an
arbitrary subset of $\{1, \dots, \ell\}$. First we recover the classes $\pm
K_0$
as the two elements of maximal square in this collection. The
$E_i$ are then the elements of the collection orthogonal to $K_0$ of square
$-1$.
Thus $\pm K_0$ is preserved by every orientation-preserving
self-diffeomorphism of $\tilde X$, and evry such diffeomorphism induces a
permutation of the set $\{\pm E_1, \dots, \pm E_\ell\}$.
Similar results hold for orientation-preserving diffeomorphisms between two
surfaces. This establishes Theorem 1.1 in this case.

Note that, if we had only kept track of the line bundles $L$ in the basic
classes and not the Spin${}^c$ structures, we would only be able to prove
Theorem
1.1 modulo 2-torsion in case
$\tilde X$ is not minimal. This is because, if $2T = 0$ in $H^2(\tilde X;
\Zee)$, then the sets $\pm K_0 \pm E$ and $\pm (K_0 + T) \pm (E+T)$ are equal.

Now suppose that $N$ is a closed negative definite 4-manifold
such that $\tilde X$ is diffeomorphic to $M\#N$ for some $M$. As we have
seen, $H_1(N; \Zee) = 0$. Let
$\{n_1,
\dots, n_k\}$ be a basis for $H^2(N; \Zee)$ such that $n_i^2 = -1$ for all $i$
and $n_i\cdot n_j = 0$ if $i\neq j$. The blowup formula for basic classes
\cite{10} implies that $M$ is of simple type and that the basic classes for
$X$ are exactly of the form
$K +\sum _{i=1}^k\pm n_i$, where $K$ is a basic class for $M$. There is also
a formula for the corresponding Spin${}^c$ structures which is analogous to
(2.5)(ii). Thus
$n_i$ is of the form
$\dsize\frac{\tilde L_1 -
\tilde L_2}2$ for two distinct basic classes $\tilde L_1$ and $\tilde L_2$, and
$n_i\in \{\pm K_0, \pm K_0 + \sum _{i\in A'} \pm E_i, \sum _{i\in A}\pm E_i\}$,
where $A'$ denotes  a proper subset of
$\{1, \dots, \ell\}$, and  $A$ denotes  an
arbitrary subset of $\{1, \dots, \ell\}$. Since $n_i^2 =-1$, the only
possibility
is $n_i = \pm E_j$ or $n_i = \pm K_0 + \sum _{i\in A}\pm E_i$,
where $\#(A) = K_0 ^2  +1$. On the other hand, the isometry given by reflection
in
$n_i$ fixes the set of basic classes $(L, \xi)$ and thus the set of differences
of basic classes. Thus the reflection must send $K_0$ to $\pm K_0$. We may
assume that $n_i = K_0 + \sum _{i\in A}\pm E_i$. Then reflection in $n_i$
applied to $K_0$ gives
$$K_0 + 2(K_0)^2(K_0 + \sum _{i\in A}\pm E_i).$$
Since $K_0^2 >0$, $K_0 + 2(K_0)^2(K_0 + \sum _{i\in A}\pm E_i)= (1+
2(K_0)^2)K_0 +  2(K_0)^2(\sum _{i\in A}\pm E_i)$ is never a difference of basic
class. This is a contradiction. Thus we must have
$n_i = \pm E_j$ for some $j$, proving Theorem 1.2 in this case.

Now let us extend the argument to handle the case $X$ is minimal and
$p_g(X) \neq 0$, but where
$K_X^2 = 0$. We again let $\tilde X$ be a blowup of $X$ with exceptional
classes
$E_i$ and denote by $K_0$ the image of $K_X$ in $H^2(\tilde X; \Zee)$. The
above
argument tells us how to recover the basic classes for
$X$ from the basic classes for $\tilde X$: they appear as the
differences $\dsize \frac{\tilde L_1-  \tilde L_2}2$ of square zero, where the
$\tilde L_i$ range over the basic classes of $\tilde X$. Moreover
the classes
$\pm K_0$ are characterized among such classes by noting that, in rational
cohomology, all other classes $K$ can be written as $rK_0$ with $|r|< 1$, at
least if $K_0$ is not torsion, whereas if $K_0$ is torsion then the only basic
classes are $\pm K_0$. Thus we see that
$\pm K_0$ is preserved by every orientation preserving
self-diffeomorphism of
$\tilde X$, and similarly for diffeomorphisms between two surfaces.

We must now recover the classes of the exceptional curves, which again appear
as differences of basic classes, and are of square $-1$ orthogonal to $K_0$.
However in this case we have additional difference classes
$\dsize\frac{K_1 - K_2}2 \pm E_i$ of square
$-1$, which are also orthogonal to $K_0$. Moreover every difference class of
square $-1$ is of the form $\dsize\frac{K_1 - K_2}2 \pm E_i$, where $K_i \in
H^2(X; \Zee)$.  Suppose that some cohomology class $\alpha$ is represented by
an
embedded
$2$-sphere, or more generally lies in
$H^2(N; \Zee)$, where $N$ is a negative definite 4-manifold such that $\tilde
X$ is diffeomorphic to $M\#N$ for some $M$. The blowup formula says that the
basic classes for $\tilde X$ are exactly of the form $\pm \alpha  \pm L$ for
certain classes $L$ orthogonal to $\alpha$. In particular $\alpha$ is again an
difference class, and it has square $-1$. Thus $\alpha =
\dsize\frac{K_1 - K_2}2 \pm E_i$ for some $K_1, K_2 \in H^2(X; \Zee)$. Set $T =
\dsize\frac{K_1 - K_2}2$. After renumbering and sign
change we may assume that $\alpha = T+E_1$. Let
$K$ be an arbitrary basic class for $X$. Since $K+E_1 + \sum _{i>1}E_i$ is a
basic class for $\tilde X$, either $K+E_1 + \sum _{i>1}E_i = T+E_1+L$ or $K+E_1
+ \sum _{i>1}E_i = -(T+E_1)+L$ for some class $L$. We claim that we cannot have
$K+E_1 +
\sum _{i>1}E_i = -(T+E_1)+L$, for otherwise we would have
$$L = K+T+2E_1 + \sum _{i>1}E_i$$
and in this case
$$T+E_1+L =  K+2T+3E_1 + \sum _{i>1}E_i$$
would be a basic class for $\tilde X$, which is clearly impossible.
Thus
$L= K-T + \sum _{i>1}E_i$. Since $ -(T+E_1)+L$ is also a basic class, we see
that
$$K-2T -E_1 + \sum _{i>1}E_i$$
is a basic class for $\tilde X$ whenever $K$ is a basic class for $X$. Now the
basic classes for $\tilde X$ are of the form $K' + \sum _i\pm E_i$, for $K'$ a
basic class on $X$. Since $T\in H^2(X; \Zee)$, the only such class
that can equal
$K-2T -E_1 +
\sum _{i>1}E_i$ is $K' -E_1 + \sum _{i>1}E_i$ for some basic class
$K'$ on
$X$. It follows that, if
$K$ is a basic class for
$X$, then  $K-2T$ is also a basic class for $X$. Since there are
only finitely many basic classes, $T$ is torsion. Now applying the above to
$K=K_0$, we see that $K_0 -2T$ is a basic class which equals $K_0$ in rational
homology. By the last sentence in (2.3), it follows that
$2T = 0$. Finally, keeping track of the Spin${}^c$ structures in the blowup
formula for taking connected sum with a negative definite 4-manifold, one
checks that the difference class
$$\frac{(L+\alpha)-(-L+ \alpha)}2 = L = K_0-T + \sum _{i>1}E_i.$$
Since $T$ is torsion, it follows from (2.3) that the only way that such a class
can be of the form
$$\frac{K_1-K_2}2 + \sum _{i\in A} \pm E_i,$$
where $K_1$ and $K_2$ are basic classes on $X$, is if $K_1=K_0$ and
$K_2=-K_0$. In this case the difference  $\dsize\frac{K_0-(-K_0)}2$ is equal to
$K_0$, and so $T=0$.

In general, using somewhat different methods, one can show the following (cf\.
\cite{15}, Chapter VI, Theorem 5.3 for the corresponding result in Donaldson
theory):

\proposition{3.1} Let $M$ be a closed oriented $4$-manifold with $b_2^+(M)
\geq 3$ and such that the set of basic classes for $M$ is nonempty. Suppose
that
$M$ is orientation-preserving diffeomorphic to $M_1
\#N_1$ and also to $M_2 \#N_2$, where the $N_i$ are negative
definite $4$-manifolds with $H_1(N_1) = H_1(N_2) = 0$. Let $n_1, \dots, n_r$ be
the exceptional classes for
$N_1$ and $n_1', \dots, n_s'$ the exceptional classes for $N_2$. Then, for
every
$i$, $1\leq i\leq s$, either there exists a $j$, $1\leq j \leq r$, such
that $n_i' = \pm n_j$ mod torsion, or $n_i'$ is orthogonal to the span of the
$n_j$. \qed
\endproclaim

However, without more knowledge about the nature of the basic classes as in the
case of a K\"ahler surface, the equality mod torsion in (3.1) seems to be the
optimal statement.

Finally let us show that the basic classes determine the plurigenera of $\tilde
X$. In all cases the basic classes determine $K_0^2$. If $K_0^2 > 0$, then
$\tilde X$ is of general type. It is well-known that, for $n\geq 2$,
$$P_n(\tilde X) = \frac{n(n-1)}2K_0^2 + \chi (\scrO_{\tilde X}).$$
Thus $P_n(\tilde X)$ is determined from the knowledge of $K_0^2$ for all $n\geq
2$, and $P_1(\tilde X) = p_g(\tilde X)$ is an oriented homotopy invariant
since $b_2^+(\tilde X) = 2p_g(\tilde X)+1$. Hence the plurigenera are
determined by $K_0^2$ as long as $K_0^2 >0$.

For
$K_0^2 = 0$,
$X$ is deformation equivalent to an elliptic surface and the plurigenera are
essentially determined from the knowledge of the multiple fibers (see
\cite{15}, Chapter I Proposition 3.22 for a more complete discussion). We will
deal here with the simply connected case, the case of at most two multiple
fibers of relatively prime multiplicity. Here the smooth classification of
elliptic surfaces with
$p_g \geq 1$ has been worked out in Donaldson theory; see
\cite{2}, \cite{23}, \cite{22}, \cite{12}, as well as \cite{21} and \cite{11}.
The general case may be reduced to the simply connected case (in the case of
finite cyclic fundamental group) or dealt with either by elementary arguments
involving the fundamental group (in case the fundamental group is not finite
cyclic), see
\cite{15} for this reduction. One could also use the basic classes to determine
the multiplicities in the general case along the lines of what we do here for
the
simply connected case.

Suppose that there are two multiple fibers $F_1$ and $F_2$, with relatively
prime multiplicities $m_1\leq m_2$. Here to handle all cases at once we will
also allow $m_1$ or both $m_1$ and $m_2$ to be 1. All the basic classes are of
the form $r\kappa$, where $\kappa$ is a primitive integral class and $r\in
\Bbb Q$. The largest value of $|r|$ is attained for $\pm K_X$, and it is
$(p_g+1)m_1m_2 -m_1 -m_2$. The next largest  value is attained
for $L = \pm (K_X -2F_2)$, and it is
$$(p_g+1)m_1m_2 -m_1 -m_2- 2m_1,$$
provided that this number is not negative. Note that since $p_g \geq 1$, if
$m_1 \geq 2$, in other words if there are two multiple fibers, then as $m_1<
m_2$,
$$(p_g+1)m_1m_2 -m_1 -m_2- 2m_1 > 4m_2 - 4m_2 = 0.$$
Thus we determine $(p_g+1)m_1m_2 -m_1 -m_2$ and $m_1$, and so
$$((p_g+1)m_1-1)(m_2 -1) = (p_g+1)m_1m_2 -m_1 -m_2 -p_gm_1 + 1.$$
{}From the knowledge of $m_1$ and $((p_g+1)m_1-1)(m_2 -1)$ we may then
determine
$m_2$.

If $m_1 = 1$ then $(p_g+1)m_1m_2 -m_1 -m_2- 2m_1 = p_gm_2 -3 \geq 0$ provided
that $p_gm_2 \geq 3$. Except in these cases, we then find $m_1 = 1$ and can
solve for $m_2$ as before. The remaining cases are $p_g = 1$, $m_2 = 1$ or $2$
or $p_g = 2$ and $m_2 = 1$. For example if $p_g = 1$, the basic classes are
$\pm
K_0$ and $K_0$ is trivial if $m_2 = 1$ and nontrivial otherwise. Thus we can
distinguish these cases also. So we have determined the multiplicities of the
multiple fibers and thus the plurigenera.

\section{4. The case where $\boldkey p_{\boldkey g}$ is zero and $\boldkey X$
is of general type.}

In this section we consider the case where $p_g(X) = 0$, i\.e\. $b_2^+(X) = 1$,
and $X$ is of general type. If $X$ is of general type, then automatically
$b_1(X) = 0$ since
$\chi (\scrO_X) = 1 - q(X) + p_g(X) > 0$.  For a $4$-manifold
$M$ with $b_2^+(M) = 1$, the function
$SW _{M,g}$ is no longer independent of the metric $g$, at least for $b_2(M)$
sufficiently large. For the purposes of this paper, we shall only consider
basic
classes of index zero. Equivalently we shall restrict the function $SW _{M,g}$
to
a function defined on characteristic line bundles $L$ with $L^2 = 2\chi (M) +
3\sigma (M)$. In this case, if
$b_2(M) = b\geq 10$,  the orthogonal hyperplanes to characteristic cohomology
classes of square
$10-b$ divide
$$\Bbb H(M) = \{\, \alpha \in H^2(M; \Ar): \alpha ^2 = 1\,\}$$
into a set of chambers $\Cal C$, and we can define $SW _{M,\Cal C}$, for each
chamber $\Cal C$, as a function on pairs of characteristic line bundles and
Spin${}^c$ structures
$(L, \xi)$ with
$L^2 = 2\chi (M) + 3\sigma (M)$. Here given a metric
$g$, there is a unique associated self-dual harmonic 2-form $\eta$ for $g$, mod
nonzero scalars. For a generic metric $g$, $t\eta$ will lie in the interior of
a
chamber $\Cal C$ for an appropriate positive real number $t$, and we let $SW
_{M,\Cal C} = SW _{M,g}$ for every metric $g$ whose self-dual 2-form lies in
$\Ar^+\cdot \Cal C$. Changing
$\eta$ to $-\eta$ corresponds to changing the orientation of the SW moduli
space.
Thus
$$SW _{M, -\Cal C} = -SW _{M, \Cal C}.$$

This procedure defines $SW _{M, \Cal C}$ for every $\Cal C$ which contains
the self-dual harmonic 2-form of a Riemannian metric. However, we can define
$SW _{M, \Cal C}$ in general formally once we have a wall crossing formula
\cite{10}:

\proposition{4.1} Suppose $b_1(M) = 0$ and that $\Cal C_0$ and $\Cal C_1$ are
two
chambers whose boundaries intersect in an open subset of a wall $L^\perp$.
Suppose further that there is a path of metrics $\{\,g_t : t \in [0,1]\,\}$
such
that the self-dual harmonic $2$-form associated to $g_0$ lies in $\Cal C_0$ and
that the self-dual harmonic $2$-form associated to $g_1$ lies in $\Cal C_1$.
Then
$$SW _{M, \Cal C_1}(L',\xi) = \cases SW _{M, \Cal C_0}(L',\xi), &\text{if
$c_1(L')\neq \pm c_1(L)$ in $H^2(M; \Bbb R)$;}\\
SW _{M, \Cal C_0}(L', \xi) \pm 1, &\text{if $c_1(L')= \pm c_1(L)$ in $H^2(M;
\Bbb
R)$.}
\endcases$$
\endproclaim

Here the sign in the wall crossing formula depends on whether  $L' \cdot \Cal
C_0$ is positive or negative, as well as on the general conventions we have
used
to orient the moduli space.  Although we shall not need to know the sign
precisely, let us give the correct choice of sign in the case of interest to
us:

\claim{4.2} Suppose in the above situation that $X$ is a K\"ahler surface and
that our orientations are chosen so that the SW moduli space has its natural
complex orientation. Suppose further that $L$ is a wall of $\Cal C_+$ and $\Cal
C_-$ and that $L\cdot \Cal C_+ >0> L\cdot \Cal C_-$. Then
$$SW _{X, \Cal C_+}(L, \xi) = SW _{X, \Cal C_-}(L, \xi) -1.$$
\endproclaim

Let us show that the claim holds in a special situation. Suppose that $X$ is a
rational surface which is the blowup of $\Pee ^2$ at $d^2$ points which lie on
a smooth curve $C$ of degree $d \geq 4$. Thus $g(C) \geq 3$. We continue to
denote by
$C$ the proper transform of $C$ on $X$. On $X$, $C^2 = 0$ and so $K_X\cdot C +
C^2 = K_X\cdot C = 2g(C) -2 >0$. Let $H$ denote the pullback of the positive
generator of $H^2(\Pee ^2; \Zee)$ to $X$ and let the classes of the
exceptional curves be denoted by $E_1, \dots, E_{d^2}$. Let
$g$ be a K\"ahler metric on
$X$ with K\"ahler form $\omega _0$ equal to $NH - \sum _iE_i$, for $N\gg 0$.
Thus
as
$K_X = -3H + \sum _iE_i$, if $N\gg0$ $\omega _0 \cdot K_X <0$. It follows that
there are no holomorphic line bundles $L$ on $X$ with $L\cdot \omega _0\leq 0$
and
$\dsize \frac{K_X+L}2 $ effective (we would simultaneously have $(K_X+L) \cdot
\omega _0 <0$ and $(K_X+L)\cdot \omega _0 \geq 0$). Thus there are no basic
classes for the chamber $\Cal C_+$ containing $\omega _0$. Since $C^2=0$, the
curve $C$ has nonnegative intersection with every irreducible curve, so by the
Nakai-Moishezon criterion
$\omega = \omega _0 + tC$ is ample for every $t$. Using the fact that
$K_X\cdot C >0$, we can choose $t$ so large that $\omega \cdot K_X >0$. It
follows that $\pm K_X$ are basic classes for the chamber $\Cal C_-$ containing
$\omega$, and moreover, by (2.4),
$SW _{X, \Cal C_-}(-K_X) = 1$ for the natural Spin${}^c$ structure and the
natural choice of complex orientation. Now
$\Cal C_-$ and
$\Cal C_+$ are separated by the wall $(-K_X)^\perp$, and $-K_X\cdot \Cal C_+ >
0
> -K_X\cdot
\Cal C_-$. Thus  $SW _{X, \Cal
C_+}(-K_X)  = SW _{X, \Cal C_-}(-K_X) - 1$.

The sign in the general case can be computed by showing that there is a
universal sign associated to the wall-crossing formula and working it out in a
special case as above.

Using the wall crossing formula, we can define $SW _{M, \Cal C}$ for all
chambers $\Cal C$, and this definition
 agrees with $SW _{M,g}$ in case $g$ is a generic
Riemannian metric whose associated harmonic self-dual 2-form lies in $\Cal C$.
A
characteristic vector
$L$ of the appropriate square for which $SW _{M, \Cal C}(L) \neq 0$ will be
called a {\sl basic class for the chamber $\Cal C$} (recall however that we
are only considering basic classes of index zero).

We consider now the following situation: $X$ is a minimal surface of general
type with $p_g(X) = 0$ and $b_1(X) = 0$.
Let $n = K_X^2$, so that $1\leq n\leq 9$. Let $\tilde X$ be a blowup of $X$ at
$\ell$ distinct points. Let $K_0$ be the pullback of $K_X$ in $\tilde X$ and
let $E_1, \dots, E_\ell$ be the exceptional curves.

\proposition{4.3} There is a unique chamber $\Cal C_0$ containing $tK_0$ in its
interior for some positive real number $t$. The chamber $\Cal C_0$ is
invariant under reflection by the classes $E_i$. The basic classes for
$\Cal C_0$ are then $\pm K_0 + \sum _{i=1}^\ell \pm E_i$, with the natural
{\rm Spin}${}^c$ structures.
\endstatement
\proof  If $K_0$ lies on a wall, then $K_0$ is orthogonal to some
characteristic cohomology element $L = A + \sum _{i=1}^\ell m_iE_i$. Here $A\in
H^2(X; \Zee) \subseteq H^2(\tilde X; \Zee)$ is characteristic, hence nonzero,
and the $m_i$ are odd, hence nonzero. Since $K_0$ is perpendicular to $L$, we
have $A^2\leq 0$. Thus $L^2\leq -\sum _im_i^2 \leq -\ell$. But on the other
hand
$L^2 = n-\ell$, which is impossible. Thus $K_0$ lies on no wall, hence lies in
the interior of a chamber $\Cal C_0$. In this case, it follows from (2.2) and
from the blowup formula (2.5) that the basic classes are as claimed.
\endproof

\corollary{4.4} The chamber $\Cal C_0$ has the following properties:
\roster
\item"{(i)}" For all characteristic $L$ with $L^2 = K_{\tilde X}^2 = n-\ell$,
$SW _{\tilde X, \Cal C_0}(L, \xi)$ is zero or $\pm 1$.  Moreover, if $(L, \xi)$
and $(L, \xi ')$ are two basic classes for $\Cal C_0$, then $\xi = \xi '$.
\item"{(ii)}" There are $2^{\ell +1}$ $(L, \xi)$ as in \rom{(i)} such that  $SW
_{\tilde X, \Cal C_0}(L, \xi)=\pm 1$.
\item"{(iii)}" If $L_1\neq L_2$ are two classes for which $SW
_{\tilde X, \Cal C_0}(L_i, \xi _i)=\pm 1$, $i=1,2$, then
$$\fracwithdelims(){L_1+L_2}{2}^2 \leq n,$$
with equality holding for at least one such pair of $L_1$ and $L_2$.
\item"{(iv)}" For every $L_1$ and $L_2$ such that equality holds above, the
line through $L_1+L_2$ meets $\Cal C_0$.  \qed
\endroster
\endproclaim

We now claim:

\proposition{4.5}
\roster
\item"{(i)}" Suppose that  $H^2(X; \Bbb Z)$ has no $2$-torsion, or more
generally that the group of $2$-torsion elements of $H^2(X; \Bbb Z)$ is not
isomorphic to $\Bbb Z/2\Bbb Z$. If
$\Cal C$ is any chamber satisfying
\rom{(i)---(iii)} of
\rom{(4.4)}, then $\Cal C = \pm \Cal C_0$.
\item"{(ii)}" If the group of $2$-torsion elements of $H^2(X; \Bbb Z)$ is
isomorphic to $\Bbb Z/2\Bbb Z$, and $\Cal C$ is any chamber satisfying
\rom{(i)---(iii)} of
\rom{(4.4)}, then $(L, \xi)$ is a basic class for $\Cal C$ if and only if
there exists a $\xi '$ such that $(L, \xi')$ is a basic class for
$\Cal C_0$.
\endroster
\endproclaim
\proof Throughout this proof we shall identify classes in $H^2(\tilde X; \Zee)$
with their images in $H^2(\tilde X; \Ar)$, i\.e\. with their images mod
torsion, since we are only concerned with the chamber structure which lives
inside $H^2(\tilde X; \Ar)$. Let
$\Cal C$ be any chamber satisfying (i)---(iii) of (4.4). Since
$-\Cal C$ also satisfies these conditions, we can assume that $\Cal C$ and
$\Cal
C_0$ lie in the same component of $\Bbb H(\tilde X)$. Now let us show that we
may
assume that $K_0 + \sum _iE_i$ is positive on both $\Cal C_0$ and $\Cal C$:

\lemma{4.6} Suppose that $\Cal C$ is a chamber lying on the same component of
$\Bbb H(\tilde X)$ as $\Cal C_0$. Then, possibly  after replacing $\Cal C$ by
$f^*\Cal C$, where
$f$ is an orientation-preserving self-diffeomorphism of $\tilde X$
corresponding
to reflection about some of the $E_i$'s, we may assume that $K_0 + \sum _iE_i$
is
positive on both $\Cal C_0$ and $\Cal C$. Thus $SW _{\tilde X,\Cal C}$ is not
identically zero for any chamber $\Cal C$. In particular $\tilde X$ is not
diffeomorphic to a rational surface and does not have a Riemannian metric of
positive scalar curvature.
\endproclaim
\proof Suppose that $x\in \Cal C$. Write $x = B + \sum _it_iE_i$, where $B\in
H^2(X;\Ar)$. After composing by a sequence of reflections in the $E_i$ (which
leave the chamber $\Cal C_0$ invariant) we may assume that $t_i \leq 0$ for all
$i$. Moreover
$x$ and
$K_0$ lie in the same component of the positive cone of $H^2(X; \Ar)$. Thus
$K_0\cdot x > 0$.  It follows that $(K_0 + \sum _iE_i)\cdot x = K_0 \cdot x +
\sum _it_i >0$. Thus  $K_0 + \sum _iE_i$ is
positive on $\Cal C$. By the wall-crossing formula, since $K_0 + \sum _iE_i$
does not separate $\Cal C_0$ from $\Cal C$, $SW_{\tilde X,\Cal C}(K_0 + \sum
_iE_i)\neq 0$.  Hence $SW _{\tilde X,\Cal C}$ is not zero on any chamber. The
final statements are then clear.
\endproof

We note that the fact that no general type surface $\tilde X$ can be
diffeomorphic to a rational surface was proved via Donaldson theory in
\cite{16} (see also \cite{19}, \cite{27}, \cite{28}, \cite{31},
\cite{32}, \cite{30}, \cite{29}). Okonek and Teleman \cite{25} have
independently observed that one can use the method of (4.6) to show that no
surface of general type is diffeomorphic to a rational surface.

Returning to the situation where $\Cal C$ satisfies (i)---(iii) of (4.4), we
see
that we can assume that $SW _{\tilde X, \Cal C}(K_0 + \sum _iE_i) \neq 0$,
and that there exists an $x\in \Cal C$ of the form $B - \sum _ir_iE_i$, where
$B\in H^2(X;\Ar)$ and $r_i \geq 0$. Now consider the walls that
separate $\Cal C$ from $\Cal C_0$. Suppose that $L= C +
\sum _i(2s_i+1)E_i$ is such a wall, where $C\in H^2(X; \Zee)$, and suppose that
$SW _{\tilde X, \Cal C_0}(L, \xi) = 0$ for all choices of $\xi$. Thus $L$ is a
basic class for
$\Cal C$ but not for $\Cal C_0$. In this case we shall show that, possibly
after
modifying $L$, the condition (iii) of (4.4) does not hold, in other words
there is an average of classes for $\Cal C_0$ with square greater than $n$.

Since
$L$ is just defined up to sign, we may assume that $L\cdot K_0 >0$, and so
$L\cdot
\Cal C_0 > 0$. Thus
$L\cdot x<0$ for all
$x\in \Cal C$. We claim that we can replace $L = C+ \sum _i(2s_i+1)E_i$ by the
class $L' = C -\sum _i|2s_i+1|E_i$. Indeed,
$$x\cdot L' = (B\cdot C) -\sum _ir_i|2s_i+1| \leq (B\cdot C) -\sum
_ir_i(2s_i+1) = x\cdot L <0< K_0\cdot L = K_0\cdot L'.$$
Hence the class $L'$ also defines a wall which separates $\Cal C_0$ and $\Cal
C$. Replacing $L$ by $L'$,  we can assume that $L= C- \sum
_i(2s_i+1)E_i$ where
$2s_i + 1 \geq 0$ and so $s_i \geq 0$. Moreover $L$ is a basic class for $\Cal
C$. It follows that $L^2 = n-\ell$ and so $C^2 = n-\ell + \sum _i(2s_i+1)^2$.
Since
$2s_i+1 \geq 1$, with equality only if $s_i=0$, we see that
$C^2 \geq n$,
with equality only if $s_i = 0$ for all $i$.

By assumption on the chamber $\Cal C$,
$((K_0 + \sum _iE_i)+ L)/2$ has square at most
$n$. If we calculate the square, however, we find:
$$\fracwithdelims(){(K_0 + \sum _iE_i)+ L}{2}^2 = \left(\frac{K_0 + C}2- \sum
_is_iE_i\right)^2= \frac{K_0^2 + C^2 +2(K_0\cdot C)}4 - \sum _is_i^2.$$
Now $K_0^2 = n$ and $C^2 = L^2 + \sum _i(2s_i+1)^2= n-\ell + \sum
_i(2s_i+1)^2\geq n$. By  the Hodge index theorem,
$$|(K_0\cdot C)| = (K_0 \cdot C) \geq \sqrt{K_0^2}\sqrt{C^2}\geq
\sqrt{n}\sqrt{n},$$
with equality holding only if $K_0 = C$. Hence
$$\align
\frac{K_0^2 + C^2 +2(K_0\cdot C)}4 - \sum _is_i^2 &\geq \frac{n + n-\ell +
\sum _i(2s_i+1)^2 +2n}4 - \sum _is_i^2\\
&= \frac14(4n + 4\sum _is_i^2 + 4\sum _is_i) - \sum _is_i^2 \geq n +\sum _is_i.
\endalign$$
Thus the square of $((K_0 + \sum _iE_i)+ L)/2$ is greater than $n$ unless $s_i
=
0$ for all $i$ and $K_0 = C$. In this case $L= K_0 - \sum _iE_i$, which is
already a basic class for $\Cal C_0$. This contradicts the choice of $L$.

It follows that the only walls which can separate $\Cal C_0$ from $\Cal C$
are are of the form $L^\perp$, where $(L, \xi)$ is a basic class for $\Cal
C_0$.
First suppose that $H^2(X; \Bbb Z)$ has no $2$-torsion. In this case we shall
show that (ii) of (4.4) does not hold, in other words that $\ell$ decreases, if
the set of such walls is nonempty. In any case the basic classes for
$\Cal C$ must be a subset of the basic classes for $\Cal C_0$. The wall
crossing
formula then implies that, if there is such a wall, then there are fewer
classes
for
$\Cal C$ than for $\Cal C_0$. This contradicts our assumption on $\Cal C$.
(Note
that, if we did not know the sign in the wall crossing formula, we would still
be
able to conclude at this point that either $\Cal C$ had fewer classes than
$\Cal C_0$ or that there existed a basic class $L$ for $\Cal C$ for which
$SW _{\tilde X, \Cal C}(L) = \pm 2$. This would again contradict the
choice of $\Cal C$.)

If there is $2$-torsion in $H^2(X; \Bbb Z)$, then the wall crossing formula
implies that we lose the basic classes $(\pm L, \xi)$ as we cross the wall
$L^\perp$ but we gain new classes of the form $(\pm L, \xi')$ for $\xi' \neq
\xi$ (recall that for $\Cal C_0$ there is a unique $\xi$ such that $(L, \xi)$
is a basic class). If the group of $2$-torsion elements of $H^2(X; \Bbb Z)$ is
larger than $\Bbb Z/2\Bbb Z$, then there would be two distinct Spin${}^c$
structures $\xi _1\neq \xi _2$ such that $(L, \xi _1)$ and $(L, \xi _2)$ are
basic classes for $\Cal C$. Thus $\Cal C$ would violate (i) of (4.4). In the
remaining case it is clear that the basic classes for $\Cal C$ are exactly of
the form $(L, \xi')$, where $(L, \xi)$ is a basic class for $\Cal C_0$.
\endproof

Note that we do not as yet claim that the integer $n$ or equivalently $\ell$ is
specified by the 4-manifold $\tilde X$, in the sense that there might {\it a
priori\/} exist other positive integers $n'\leq 9$ and $\ell'$ with $n' - \ell
' = n-\ell$, and a chamber $\Cal C_0 '$ satisfying (i)---(iii) of (4.4) for
$n'$ and $\ell'$. In fact, although we shall not really need this, our
arguments
show that $\ell$ is the maximum over all possible $\ell '$ such that there
exists
a chamber $\Cal C_0 '$ satisfying  (i)---(iii) of (4.4) for $\ell'$ and $n' =
K_{\tilde X}^2 + \ell$.

Let us now show how to deduce Theorems 1.1 and 1.2 from Proposition 4.5, at
least in the case where $p_g(X) = 0$ and $X$ is of general type. First suppose
that there is no $2$-torsion in $H^2(X; \Bbb Z)$ and that $f\: \tilde X \to
\tilde X$ is an orientation-preserving self-diffeomorphism. Then $f^*\Cal C_0$
has the same properties (i)---(iii) as
$\Cal C_0$, and so by (4.5) $f^*\Cal C_0 = \pm \Cal C_0$. In particular $f^*$
leaves invariant the basic classes for $\Cal C_0$. As in the case $p_g>0$ and
$K_0^2>0$ we see that $f^*$ preserves $\pm K_0$ and the span of the $E_i$. A
similar argument works in case the $2$-torsion subgroup of  $H^2(X; \Bbb Z)$
is not isomorphic to $\Bbb Z/2\Bbb Z$. Finally, if the $2$-torsion subgroup of
$H^2(X; \Bbb Z)$ is isomorphic to $\Bbb Z/2\Bbb Z$ and $f\: \tilde X \to
\tilde X$ is an orientation-preserving self-diffeomorphism, then $f^*\Cal C_0$
has the same properties (i)---(iii) as
$\Cal C_0$, and so  $f^*$
leaves invariant the set of $L$ such that there exists a $\xi$ with $(L,
\xi)$ a basic class for
$\Cal C_0$. Thus $f^*$ preserves $\pm 2K_0$ and so $f^*\Cal C_0 = \pm \Cal
C_0$. Again we may complete the argument as in the case $p_g >0$.

Next, let $N$ be a  negative definite 4-manifold
such that $\tilde X$ is diffeomorphic to $M\#N$ for some 4-manifold $M$. Let
$x\in H^2(M; \Ar) \subseteq H^2(\tilde X; \Ar)$ satsify $x^2 =1$, and assume
that $x$ lies on no wall in $\Bbb H(M)$ or $\Bbb H(\tilde X)$. Let $\Cal C$ be
the chamber of $\Bbb H(\tilde X)$ containing $x$. Let $n_i$ be an exceptional
class of $N$ and let $R_i\: H^2(\tilde X; \Ar) \to H^2(\tilde X;
\Ar)$ be the reflection about $n_i$. Note that, as $H^2(N; \Zee)$ is torsion
free, $R_i$ acts on the set of Spin${}^c$ structures as well. Clearly
$R_i$ fixes
$\Cal C$. The blowup formula says that $R_i$ also fixes the set of basic
classes
for
$\Cal C$.

We claim that, in the above situation, the basic classes for $R_i(\Cal C_0)$
are
exactly those of the form $R_i(L)$, where $L$ is a basic class for $\Cal C_0$,
and more precisely that $SW _{\tilde X, R_i(\Cal C_0)}(R_i(L), R_i(\xi)) = SW
_{\tilde X,\Cal C_0}(L, \xi)$. (This would of course be clear if $n_i$ was
represented by a smoothly embedded 2-sphere.) Thus, arguing as above,
$R_i(\Cal C_0)$ satisfies (i)---(iii) of (4.4) and moreover
$R_i(\Cal C_0) =
\pm \Cal C_0$, indeed $R_i(\Cal C_0) = \Cal C_0$ since $R_i$ fixes the
components of $\Bbb H(\tilde X)$. It follows that the wall through $n_i$ passes
through $\Cal C_0$ for every $i$, and that (by induction on the number of
exceptional classes) we can choose an
$x\in
\Bbb H(M)$ whose image in  $\Bbb H(\tilde X)$ lies in $\Cal C_0$. The proof
that every exceptional class $n_i$ must be equal to $\pm E_j$
for some $j$ then runs as in the case
$p_g(X)>0$ and $X$ of general type (i\.e\. $K_0^2 >0$). Thus it suffices to
show (we omit the Spin${}^c$ structures for notational simplicity):

\proposition{4.7} Let $Y=M\#N$ be an oriented smooth $4$-manifold such that
$b_2^+(M) =1$ and $N$ is negative definite. For an exceptional
class $n_i \in H^2(N; \Zee)$, let $R_i$ be the refection about $n_i$, viewed as
an automorphism of $H^2(Y; \Zee)$. Then for every chamber $\Cal C_0$ of $\Bbb
H(Y)$, we have:
$$SW _{Y, R_i(\Cal C_0)}(R_i(L)) = SW
_{Y,\Cal C_0}(L).$$
\endproclaim
\proof As above, let
$x\in H^2(M; \Ar) \subseteq H^2(Y; \Ar)$ satsify $x^2 =1$, and assume
that $x$ lies on no wall in $\Bbb H(M)$ or $\Bbb H(Y)$. Let $\Cal C$ be
the chamber of $\Bbb H(Y)$ containing $x$, so that the blowup formula holds for
the basic classes of $\Cal C$. To see that $SW _{Y, R_i(\Cal C_0)}(R_i(L)) = SW
_{Y,\Cal C_0}(L)$, we
may assume that $\Cal C$ and $\Cal C_0$ lie on the same component of $\Bbb
H(Y)$.  Choose an $x_0 \in \Cal C_0$ and consider a generic path
$\gamma$ from
$x$ to $x_0$. Let $L_1, \dots, L_k$ be the walls crossing $\gamma$. Then if
$L\neq L_i$ for any $i$, the wall-crossing formula says that
$$SW _{Y,\Cal C_0}(L) =SW _{Y,\Cal C}(L).$$
If $L=L_\alpha$ for some $\alpha$, then
$$SW _{Y,\Cal C_0}(L_\alpha) =SW _{Y,\Cal
C}(L_\alpha)\pm 1.$$ Now consider the function $SW _{Y, R_i(\Cal
C_0)}$. The point
$R_i(x_0)$ lies in $R_i(\Cal C_0)$, the path $R_i(\gamma)$ joins $R_i(x)$ to
$R_i(x_0)$, and the walls crossed by $R_i(\gamma)$ are the walls
$R_i(L_\alpha)$. Note that the blowup formula implies that $L$ is a basic class
for $\Cal C$ if and only if $R_i(L)$ is a class for $\Cal C$, and indeed
$SW _{Y, \Cal C}(L) = SW _{Y, \Cal C}(R_i(L))$. Now to
calculate $SW _{Y, R_i(\Cal C_0)}(R_i(L))$, first assume that
$L\neq L_\alpha$ for any $\alpha$. Then
$$SW _{Y, R_i(\Cal C_0)}(R_i(L)) = SW _{Y, \Cal
C}(R_i(L)) = SW _{Y, \Cal C}(L) = SW _{Y,\Cal C_0}(L).$$
If $L= L_\alpha$, then
$$SW _{Y, R_i(\Cal C_0)}(R_i(L_\alpha)) = SW _{Y, \Cal
C}(R_i(L_\alpha)) \pm 1 = SW _{Y, \Cal C}(L_\alpha)\pm 1.$$
Now examination of the wall crossing formula says that
$$SW _{Y,\Cal C_0}(L_\alpha) = SW _{Y, \Cal
C}(L_\alpha)\pm 1,$$
and the sign must be the same as in the above formula. Putting this together we
see that $SW _{Y, R_i(\Cal C_0)}(R_i(L_\alpha)) = SW
_{Y,\Cal C_0}(L_\alpha)$, so the formula holds in all cases.
\endproof

Thus we have established Theorem 1.2. In particular, the cohomology
classes of embedded 2-spheres of self-intersection $-1$ span a sublattice of
$H^2(\tilde X; \Zee)$ of rank exactly $\ell$. It follows that if $\tilde{X'}$
is another  surface of general type and $f\: \tilde{X'} \to
\tilde X$ is an orientation-preserving diffeomorphism, then $\tilde {X'}$ can
be
blown up at most
$\ell$ times from its minimal model. By symmetry
$\tilde X$ and $\tilde{X'}$ are blown up the same number of times, namely
$\ell$, from their minimal models. If $\Cal C_0'$ is the chamber on
$\tilde{X'}$
corresponding to $\Cal C_0$, it then follows that $f^*\Cal C_0' = \pm \Cal
C_0$. Theorem 1.1 is an immediate consequence.
\qed
\medskip

Finally let us deduce that if $Y$ is a complex surface diffeomorphic to $\tilde
X$, then $P_n(Y) = P_n(\tilde X)$ for all $n$. Note that $Y$ is K\"ahler since
$b_1(Y) =0$. Moreover $Y$ cannot be  a rational surface by  Lemma 4.6. We
will rule out the case where $Y$ is elliptic in the next section. Thus
$Y$ is again a surface of general type, and by Theorem 1.1 we may determine
$K_0^2$ for $Y$. It follows as in the case where $p_g>0$ that $P_n(Y) =
P_n(\tilde X)$ for all $n$.

We note that Theorem 1.1 has the following corollary:

\corollary{4.8} Let $\tilde X$ be a surface of general type with $p_g(\tilde X)
= 0$, and let
$D(\tilde X)$ be the image of the group of orientation-preserving
diffeomorphisms  of
$\tilde X$ in the automorphism group of $H^2(\tilde X; \Zee)$. Then
$D(\tilde X)$ is finite.
\endproclaim
\proof Let $\phi \in D(\tilde X)$. Up to finite index we may assume that
$\phi(E_i) = E_i$ for all $i$ and that $\phi (K_0) = K_0$. Thus $\phi$ is
determined by its action on $K_0^\perp \cap H^2(X; \Zee)$. Since $K_0^\perp
\cap
H^2(X; \Zee)$ is negative definite, there are only finitely many automorphisms
of $K_0^\perp \cap H^2(X; \Zee)$. Hence there are only finitely many
possibilities for $\phi$.
\endproof

\section{5. The case where $\boldkey p_{\boldkey g}$ is zero and $\boldkey X$
is
not of general type.}

Suppose that $X$ is a minimal K\"ahler surface, not of general type, rational,
or ruled, such that
$p_g(X) = 0$. In this case $X$ is elliptic (possibly an Enriques or
hyperelliptic surface), $K_X^2=0$,  and, since $\chi (\scrO_X) \geq 0$, either
$b_1(X) = 0$ and $X$ is elliptic or $b_1(X) = 2$. We shall first consider the
parts of the analysis of the topology of a blowup of $X$ which can be handled
either by elementary methods or by reduction to the case $p_g > 0$.

If $b_1(X) = 0$, then
$X$ is an elliptic surface, obtained from a rational elliptic surface by a
number of logarithmic transforms. Here $X$ is rational if there is just one
logarithmic transform and it is simply connected (a Dolgachev surface) if there
are two logarithmic transforms of relatively prime multiplicities. In all other
cases
$X$ has a finite covering space $Y$ which is an elliptic surface with $p_g \geq
1$. Let
$\tilde X$ be a blowup of $X$, and let $\tilde Y$ be the induced cover, which
is a blowup of $Y$. In this case, since Theorems 1.1 and 1.2 hold for
$\tilde Y$, they also hold for $\tilde X$ mod torsion. We can then determine
the plurigenera of $\tilde X$, either via elementary arguments involving the
fundamental group as in
\cite{15} if $\pi _1(X)$ is not finite cyclic or by reducing to the simply
connected case with $p_g \geq 1$ in case $\pi _1(X)$ is finite cyclic (but
nontrivial). Finally, since $\tilde Y$ has no Riemannian metric with positive
scalar curvature, the same is true for $\tilde X$.

Likewise, if $X$ is a
minimal K\"ahler surface with $p_g(X) = 0$ and $b_1(X) = 2$, or in other words
$q(X) = 1$, then $X$ is an elliptic surface with Euler number zero. In this
case either $X$ is ruled over an elliptic base or it is nonruled. If $X$ is not
ruled it is a logarithmic transform either of an elliptic surface without
singular or multiple fibers over a base curve of genus one, with nontrivial
holonomy (in other words, a hyperelliptic surface), or it is a logarithmic
transform of $E\times \Pee ^1$ so that the corresponding base orbifold is not
spherical (see \cite{15} for more description). In both of these cases,
provided that $X$ is not ruled, it again has a finite  covering space $Y$ which
is an elliptic  surface with $p_g \geq 1$. Again, we see that Theorems 1.1 and
1.2 hold for blowups of $X$, mod torsion. The determination of the plurigenera
then follows from the fundamental group as in \cite{15}, and the fact that
$\tilde X$ has no metric of positive scalar curvature, for every blowup $\tilde
X$ of $X$, follows as in the case where
$b_1(X) =0$ and $\pi _1(X)$ is not trivial.

The remaining case is the case of Dolgachev surfaces. This case was handled via
Donaldson theory in \cite{13}, \cite{1}, \cite{12}  (see also \cite{26}). It
can
also be handled by the methods of this paper, as we now outline. Let $\tilde X$
be the blowup of a minimal Dolgachev surface $X$, and let $K_0$ be the image of
the canonical class of $X$ in
$H^2(\tilde X; \Zee)$. Suppose that $\tilde X$ is the blowup of the minimal
surface $X$ at $\ell$ distinct points.  All of the basic classes for
$X$ are rational multiples $rK_0$ of $K_0$ with $|r| \leq 1$. We suppose that
there are
$d$ basic classes for $X$. It is easy to check that on each basic class of $X$
the value of $SW_X$, which does not depend on a choice of chamber, is $\pm 1$.

Now let $\Cal C_0$ be a chamber in $\Bbb H(\tilde X)$ which contains classes
of the form $\omega$, where $\omega$ is the  K\"ahler form of a generic
K\"ahler metric on $X$. It follows that there are $d2^\ell$ basic classes for
$\Cal C_0$. They  are exactly the classes $L + \sum _i\pm E_i$, where $L$ is a
basic class for
$X$ and the
$E_i$ are the exceptional classes on $\tilde X$, and the value of $SW_{\tilde
X,
\Cal C_0}$ on each such class is $\pm 1$. Moreover, for every average $\dsize
\frac{L_1+L_2}2$ of basic classes for  $\Cal C_0$, we have $\dsize
\fracwithdelims(){L_1+L_2}{2}^2 \leq 0$, with equality holding for some pair of
basic classes $L_1 \neq \pm L_2$. Arguments very similar to those given in the
proof of Proposition 4.5 show that
$\pm \Cal C_0$ are the unique chambers with these properties. Thus every
orientation-preserving self-diffeomorphism $f\: \tilde X \to \tilde X$
satisfies $f^*\Cal C_0 = \pm \Cal C_0$, and if $n_i$ is an exceptional class
for a negative definite summand of $\tilde X$, then the reflection $R_i$ in
$n_i$ preserves $\Cal C_0$. Moreover, in the above two cases, we see that
both $f^*$ and $R_i$ permute the set of basic classes. In particular the wall
through $n_i$ passes through $\Cal C_0$ and $n_i$ is a difference class
for $\Cal C_0$. It follows then from the arguments in Section 3 for the case
$K_0^2 =0$ that $f^*$ preserves $\pm K_0$ and that $n_i = \pm
E_j$ for some $j$. In
particular there can be at most $\ell$ disjoint smoothly embedded 2-spheres
of self-intersection $-1$ on $\tilde X$. Thus $\tilde X$ cannot be
diffeomorphic to a blown up surface of general type. The basic classes for
$\Cal C_0$ determine the the basic classes for
$X$, by arguments similar to those in Section 3 for the case where $K_X^2 =0$.
Finally, the basic classes for $X$ determine the multiplicities of the
multiple fibers, by arguments along the lines of those given in Section 3 for
simply connected elliptic surfaces with $p_g\geq 0$. Here, in case $p_g =0$,
there are a few extra cases to consider. Lastly, the arguments used to prove
Lemma 4.6 also show that there is no chamber $\Cal C$ where the function
$SW_{\tilde X, \Cal C}$ is identically zero, and thus $\tilde X$ has no metric
of positive scalar curvature. Thus we have completed the proofs of Corollary
1.4 and Corollary 1.5.

Finally we note that one can modify the proofs of the results in Section 4 to
handle the non-simply connected elliptic surfaces with $p_g =0$, and replace
equality mod torsion with equality in the non simply connected
case. The main point is to handle wall crossings in case $b_1(X) = 2$. However,
we shall not give these arguments here.

\section{6. Some open questions.}

\noindent {\bf The non-K\"ahler case.} Can one generalize the above results to
the non-K\"ahler case? The structure of non-K\"ahler complex surfaces of
nonnegative Kodaira dimension is well-understood. In particular they are all
elliptic surfaces with odd first Betti number and Euler number zero. Elementary
methods \cite{15} show that all such surfaces are $K(\pi, 1)$'s and
hence that the classes of exceptional curves are preserved under
diffeomorphisms. It is also straightforward to show that the class of a
general fiber is preserved, and so $K_X$ mod torsion, and that the plurigenera
are diffeomorphism invariants. However, it does not seem possible to extend
these
methods to handle negative definite connected summands. On the other hand, the
analysis of the SW equations for K\"ahler or symplectic manifolds may admit a
straightforward generalization to this case (possibly under some assumptions
on the metric).

\medskip
\noindent {\bf Ruled surfaces.} First we note that the Seiberg-Witten theory
accounts for all the self-diffeomorphisms of a rational surface $X$. Let
$X$ be the blowup of $\Bbb P^2$ at $\ell$ distinct points. Inside
$\Bbb H(X)$ there is a distinguished convex set, the ``super $P$-cell"
$\bold S_0$ defined in \cite{13}. Its walls are certain characteristic elements
of $H^2(X; \Zee)$ of square $K_X^2$. In fact one can show that $\bold S_0$ is
a chamber for the set of walls defined by primitive characteristic elements of
square $K_X^2$. It is shown in
\cite{13} that an automorphism $\varphi$ of the lattice $H^2(X; \Zee)$ is the
image of an orientation-preserving self-diffeomorphism of $X$ if and only if
$\varphi (\bold S_0) = \pm \bold S_0$. This result can be established by
Seiberg-Witten theory as well. It is elementary to show that every $\varphi$
such that $\varphi (\bold S_0) = \pm \bold S_0$ is the image of a
diffeomorphism, and the difficult part of the argument is to show that, for
every
orientation-preserving self-diffeomorphism $\psi$ of $X$, $\psi ^* (\bold S_0)
=
\pm \bold S_0$. Let $\Cal C_0$ be the chamber of $\Bbb H(X)$, for the set of
all
walls defined by  characteristic  elements
of $H^2(X; \Zee)$ of square $K_X^2$, which contains $\omega _0$, where $\omega
_0$ is a K\"ahler metric on $\Bbb P^2$, or equivalently contains K\"ahler
metrics on $X$ with K\"ahler form a positive multiple of $N\omega _0 - \sum
_iE_i$. Thus $\Cal C_0$ contains the period points of metrics with positive
scalar curvature, and so $SW_{X, \Cal C_0}$ is identically zero (this also
follows from the blowup formula). Moreover $\pm\Cal C_0$ are the unique
chambers $\Cal C$ such that $SW_{X, \Cal C}$ is identically zero. Hence $\psi
^*\Cal C_0 = \pm \Cal C_0$. For simplicity assume that $\psi
^*\Cal C_0 = \Cal C_0$. Now the interior of $\Cal C_0$ is nonempty and is
contained in the interior of $\bold S_0$, since the walls defining $\bold S_0$
are a subset of the walls defining $\Cal C_0$ and $\omega _0 \in \Cal C_0\cap
\bold S_0$. Thus $\psi ^*(\bold S_0)$ is a super $P$-cell such that $\psi
^*(\bold S_0) \cap \bold S_0 \neq \emptyset$. It follows by Lemma 5.3(e) on
p\. 339 of \cite{13} that
$\psi ^*(\bold S_0) = \bold S_0$.

Similar arguments imply that, if $N$ is a negative definite summand of $X$ and
$n_i$ is an exceptional class for $N$, then the reflection in $n_i$ preserves
$\bold S_0$. Using this fact and the method of proof of \cite{16}, Theorem 1.7,
one can show that every negative definite summand $N$ of $X$ can be accounted
for in the following sense: if $X$ is orientation-preserving diffeomorphic to
$M\#N$, where $N$ is negative definite, then there is an
orientation-preserving diffeomorphism $f\: X \to \tilde Y$, where $\tilde Y$
is the blowup of a rational surface $Y$, such that $H^2(N; \Zee)$ corresponds
to the span of the exceptional classes of the blowup $\tilde Y \to Y$.

If one tries to extend the above results to surfaces $X$ which are (not
necessarily minimal) ruled surfaces over a nonrational base curve, most of the
results extend with elementary proofs. For example, let $f$ be the class of a
general fiber of $X$ and $E_1, \dots, E_\ell$ be the classes of the
exceptional curves. Then the cohomology classes of smoothly embedded
2-spheres of self-intersection zero are exactly the classes $nf, n\in \Zee$,
and the cohomology classes of smoothly embedded
2-spheres of self-intersection $-1$ are  the classes $nf+E_i, n\in \Zee$.
Moreover every orientation-preserving self-diffeomorphism $\psi$ of $X$
satisfies $\psi ^*f = \pm f$, and indeed an automorphism $\varphi$ of $H^2(X;
\Zee)$ is the image of an orientation-preserving self-diffeomorphism if and
only
if $\psi ^*f = \pm f$. However, it is not clear how to generalize these results
to arbitrary negative definite summands. Such generalizations would
presumably follow from Seiberg-Witten theory provided that we have a better
understanding of the transition formula in case $b_1(X) \neq 0$. Working out
such transition formulas is of course an interesting problem in its own right.

\medskip
\noindent {\bf Questions seemingly inaccessible to Seiberg-Witten theory.} The
overall moral of the above is that the sum total of the basic classes gives us
information about the obvious invariants of a complex surface, the exceptional
curves and the pullback of the canonical class of the minimal model, and no
more. Presumably the same is true of Donaldson theory. One can ask if there is
more to the smooth topology of a complex surface than this. For example, there
exist minimal surfaces of general type, say $X_1$ and $X_2$, and an isometry
$\varphi \: H^2(X_1; \Zee) \to H^2(X_2; \Zee)$ such that $\varphi (K_{X_1}) =
K_{X_2}$, and such that $X_1$ and $X_2$ are not deformation equivalent. As
has been pointed out by Fintushel and Stern, certain pairs of Horikawa surfaces
are among the simplest examples of surfaces with this property. In particular
we
cannot distinguish such pairs
$X_1$ and
$X_2$ via the basic classes, and, at least in the case of Horikawa surfaces,
they
have the same Donaldson polynomials as well, as has been shown by the second
author and Z. Szab\'o. Can one find new smooth invariants which will show that
$X_1$ and
$X_2$ are not diffeomorphic?

In a similar vein, are there further restrictions
on self-diffeomorphisms $f$ of a K\"ahler surface $X$ of nonnegative Kodaira
dimension beyond the conditions $f^*K_0 = \pm K_0$ and $f^*E_i = \pm E_j$? For
example, if the algebraic geometry of $X$ imposes a fundamental asymmetry on
$X$, is this seen by the smooth topology?  One example of this might be a
surface which is the double cover of $\Pee ^1\times \Pee ^1$ along an
asymmetric branch divisor in the linear system $|2af_1 + 2bf_2|$, where the
$f_i$ are the fibers of the two different projections and $a,b$ are positive
integers with $a\neq b$. Is it true that for general $a$ and $b$ every
orientation-preserving self-diffeomorphism of $X$ preserves the pullbacks to
$X$ of
$f_1$ and $f_2$ up to sign (and not just $K_X$ which is a positive combination
of the pullbacks)? Even for simply connected elliptic surfaces, there is a gap
of finite index between those isometries of $H^2(X; \Zee)$ known to arise from
self-diffeomorphisms and the restrictions placed on such isometries by
Donaldson theory or Seiberg-Witten theory (see for example \cite{15}, Chapter
II, Theorem 6.5 and \cite{13}). Can one close this
gap by constructing more diffeomorphisms, or can one find new invariants which
rule out the existence of such diffeomorphisms?

\enddocument